\documentclass[review]{elsarticle}

\usepackage{hyperref}
\usepackage{indentfirst}
\usepackage{graphicx}
\usepackage{float}
\usepackage{subfig}
\usepackage{amssymb}
\usepackage{stmaryrd}
\usepackage{amsfonts}
\usepackage{amsmath}
\usepackage{latexsym}
\usepackage{color}
\usepackage{CJK}
\usepackage{mathrsfs}
\usepackage{caption}
\usepackage{algorithm}
\usepackage{algpseudocode}
\usepackage{epstopdf}
\usepackage{multirow}
\usepackage{mathptmx}
\usepackage{booktabs}

\begin{document}

\begin{frontmatter}

\title{Global dynamics for a class of discrete fractional epidemic model with reaction-diffusion}
\author[mymainaddress]{Zhenzhen Lu}


\author[mymainaddress]{Yongguang Yu}

\author[mymainaddress]{Weiyi Xu}

\author[mymainaddress]{Guojian Ren}

\author[mymainaddress]{Xiangyun Meng}
\cortext[mycorrespondingauthor]{Corresponding author}
\ead{xymeng1@bjtu.edu.cn}

\address[mymainaddress]{Department of Mathematics, Beijing Jiaotong University, Beijing, 100044, P.R.China}

\begin{abstract}
In recent years, discrete fractional epidemic models with reaction-diffusion have become increasingly popular in the literature, not only for its necessity of numerical simulation, but also for its defined physical processes. In this paper, by second order central difference scheme and L1 nonstandard finite difference scheme, a discrete counterpart of time-fractional reaction-diffusion epidemic model with generalized incidence rate is considered. More importantly, the main idea in choosing an nonstandard finite-difference scheme is to obtain unconditionally positivity in the proposed system, which leads to the proposal of the discrete epidemic model with time delay. Furthermore, the global properties of the proposed discrete system are studied, including the global boundedness of positive solutions, the existence and the global stability of equilibrium points, which are consistent with the corresponding continuous systems. Meanwhile, it shows that L1 nonstandard finite difference scheme and second order central difference scheme can keep the properties of the corresponding continuous system well. It is worth noting that, different from discrete epidemic model with integer-order, the memory Lyapunov function is constructed in this paper, which depends on the previous historical information of the proposed system. This is consistent with the non-local property of Caputo fractional derivatives. Furthermore, the $\alpha$-robustness of L1 scheme in Lyapunov function and the stability of the proposed system are guaranteed, that is the Lyapunov function constructed in this paper is consistent with the integer-order system as the fractional-order $\alpha\to 1^-$, and thus, the stability of the proposed system can be guaranteed as the fractional-order $\alpha\to 1^{-}$. Finally, numerical results are given to verify the theoretical results.
\end{abstract}

\begin{keyword}
Time-Fractional discrete-time \sep Reaction-diffusion \sep L1 nonstandard finite difference scheme \sep Time delay \sep Lyapunov function \sep Global stability \sep $\alpha$-robust

\end{keyword}

\end{frontmatter}


\section{Introduction}
Epidemiology is the investigation of why and how certain diseases occur in a particular population. Insights gained from epidemic models are important in preventing the spread of infectious disease. An earlier study on the human epidemiology was studied by Bernoulli about 250 years ago \cite{1760}. Later, a classical SIR epidemic model (S, I and R represent susceptible, infected and recovered individuals, respectively) was presented by Kermack and McKendrick, which studied the plague outbreaks in London during 1665-1666 and in Mumbai in 1906 \cite{2017W}, respectively. Boosted by the work of Kermack and McKendrick, many researchers established various epidemic models to study the spread of different infectious diseases \cite{C2010Complete,2001GLOBAL,2021The}. In particular, the dynamics of these models can be expressed as a set of coupled ordinary differential equations with integer-order \cite{1}. However, the spread of some infectious diseases depends on on both the current state and the history state, and thus cannot be adequately dealt with classic models \cite{2013Patterns}. A more recent approach has been to incorporate a history dependence into the dynamics by generalising the classic coupled ODEs using Caputo fractional-order derivatives \cite{2012Solutions,2013A}. Meanwhile, Smethurst et al. found that the waiting time between doctor visits for a patient follows a power-law model \cite{2001Are}, which can be deduced the fractional-order system. Furthermore, Angstmann et al. derived a fractional-order infectivity SIR model from a stochastic process that incorporated a time-since-infection dependence on the infectivity of individuals, and they found the fractional derivative appears in the generalised master equations of CTRW (continuous time random walk) through SIR compartments, with a power-law distribution in the infectivity \cite{2016A}. Meanwhile, Angstmann et al. showed how fractional differential operators arise naturally in these models whenever the recovery time from the disease is a power-law distributed \cite{1}. Accordingly, there has been a large number of literature considering fractional-order epidemic models \cite{2020Global,2020Analysis,2018A}. For example, a nonlinear fractional order SIR epidemic model is established along the memory with Crowley-Martin type functional response and Holling type-II treatment rate in \cite{2020Global}, a fractional hepatitis B epidemic model with the general incidence rate is considered in \cite{2020Analysis}.

Meanwhile, classic epidemic models assume a uniform population distribution in habitats, dependent only on time evolution, which means ordinary differential equations are considered in it. However, one of the fundamental challenges of epidemiology is to determine how population structure affects the spread of infectious disease, and one important aspect is spatial structure, such as influenza, measles, severe acute respiratory syndrome (SARS) and the extremely severe spread of COVID-19, which is largely driven by human travel and sustained by various means of transport. Based on the movement of people in a region, or even within a country, or even in the world, the traditional ordinary differential model is not enough to describe the spread of infectious disease realistically. Therefore, the effect of spatial diffusion should not be ignored in the study of infectious diseases. Presently, there have been numerous research works devoted to the investigation of the roles of different types of reaction-diffusion on the dynamic of infectious diseases \cite{2021Dynamical,2012Complex,2012Stability} and therein. However, in spatial diffusion papers mentioned, the standard reaction-diffusion is always considered, that is, the diffusion process of infection diseases is modeled as Brown motion, which is unreasonable. For example, the diffusion subjected to a power tail distribution is no longer the standard diffusion process, but subdiffusion or superdiffusion. One of the most important advances was proposed by Henry et al., who derived a fractional-order reaction-diffusion equation with time memory based on CTRW and verified the one-dimensional single-species fractional-order reaction-diffusion process \cite{2000Fractional}. Subsequently, a number of literatures consider fractional reaction-diffusion models of infectious diseases \cite{2022ctrw,2008Traveling,2017A,2021Dynamics,Kuniya2018Global,2021luzhenzhen}. Particularly, Zhenzhen et al. provided a physical derivation of an SIR model that includes time-fractional space-fractional, time-fractional reaction-diffusion or fractional diffusion, and they investigated how anomalous diffusion naturally introduces epidemic model \cite{2022ctrw}. Zhao et al. was concerned with the spatial propagation of an SIR epidemic model with nonlocal diffusion \cite{2021Dynamics}. Kuniya et al. proposed an SIR epidemic model with nonlocal diffusion and obtained different conditions to guarantee the global asymptotic stability of constant equilibrium points and the uniform persistence of their model \cite{Kuniya2018Global}. Zhenzhen et al. investigated the global dynamics for a class of multi-group SIR epidemic model with time fractional-order and reaction-diffusion \cite{2021luzhenzhen}.

Besides, it is impossible to obtain explicit solution for fractional-order differential equations with diffusion and only numerical solution can be studied, except in special case. Not only that, the discrete model has its own specific physical process \cite{2015from,2015discrete,2016A}. For example, Angstmann et al. considered the DTRW (discrete time random walk) to establish the discrete SIR epidemic model, and found that the diffusion limit is a continuous SIR system with the corresponding Riemann-Liouville derivative \cite{2016A}. More importantly, infectious disease statistics is not continuous but under a certain amount of time to each.  Moreover, the discrete model of infectious diseases after can be in any time step for the unit, so after discretization of epidemic model than the continuous can be used more convenient and effective statistics. There are many numerical schemes used to construct discrete systems, such as Euler scheme \cite{2020On} and nonstandard finite difference method (NSFD) \cite{2012Bifurcations}. Meanwhile, Dimitrov et al. pointed out that NSFD can keep the physical properties of the original continuous system, which is to obtain unconditional stability and positivity in the variables representing the subpopulations \cite{article}. Then Jang et al. used NSFD to construct a discrete model of infectious immigrants and proved the global stability of the equilibrium point \cite{2003Difference}. Skiguchi et al. considered a discrete SIR epidemic model with delay, and studied the global stability of the disease-free equilibrium point and the persistence of the system \cite{2010Global}. Villanueva et al. analyzed a discrete infectious model of a respiratory syncytial virus \cite{J2008Nonstandard}. However, the above studies of discrete epidemic models either did not account for any kinds of diffusion, or they only considered standard diffusion and did not account for anomalous diffusion.

To the best of our knowledge, no papers consider discrete time-fractional epidemic models with reaction-diffusion. Furthermore, L1 scheme is usually used to discretize Caputo fractional-order derivatives \cite{2020Barrier}, which can also keep the properties of Caputo fractional derivative well \cite{2021A}. The terminology `1' refers to the numerical approximation of Caputo derivatives of order $\alpha\in(0,1)$ and the `L' is perhaps by association with Riemann-Liouville \cite{1974The}. Furthermore, the L1 schemes are among the most popular numerical approximations for Caputo derivatives, and are easy to implement with acceptable precision and good numerical stability \cite{2021M}. In addition, since the L1 scheme approximates the classical first derivative by standard forward Euler difference scheme, it provides a good basis for various numerical approximations for Caputo derivatives in \cite{2019error}. Besides, L1 scheme proved that the solution of a system converges, uniformly on the space-time domain, to the solution of the classical initial-boundary value problem as $\alpha\to 1^-$, which implies L1 scheme is $\alpha$-robust \cite{2020Blow}. More important, based on a DTRW, L1 scheme can be used to describe the time evolution of particles created  by heavy tail distribution \cite{2022discretelu}. Meanwhile, NSFD is usually used to ensure that the population described by the infectious disease model is unconditionally positive \cite{Mickens2002Nonstandard}. Based on the above discussion, by L1 nonstandard finite difference scheme and second order central difference scheme, a discrete counterpart of time-fractional reaction-diffusion epidemic model with time delay is considered to make better use of real-data. Also, the proposed discrete system can be simplified to integer-order discrete system with forward difference when the fractional-order $\alpha\to 1^-$ in \cite{2018Global}. Finally, the global properties of the above discretized systems, including boundedness of positive solutions and global stability of equilibrium points, are considered in this paper. It is worth noticed that in studying the stability of two equilibrium points, Lyapunov function is strictly dependent on the historical information of the proposed system, and the stability of equilibrium point cannot be completely obtained if only relying on current information of the corresponding discrete system. This is consistent with the memorability of continuous Caputo derivatives. What's more, with the Lyapunov function and the stability of the proposed system, the $\alpha$-robustness of L1 scheme is guaranteed between the fractional-order and integer-order system, that is the Lyapunov function constructed in this paper is consistent with the integer-order as the fractional-order $\alpha\to 1^-$ \cite{2018Global}, so the stability of the proposed system can also be guaranteed with $\alpha\to 1^-$.

Based on the above  analysis, a fractional discrete SIR model for reaction-diffusion with generalized incidence rate is considered. The main contributions of this study are as follows:
\begin{itemize}
\item[$\bullet$]
Based on the monotonicity of discrete systems, the global asymptotic stability of autonomous systems with
L1 scheme is investigated by Lyapunov's direct method.
\item[$\bullet$]
The proposed discrete system is considered, in which L1 scheme is adopted in the time direction and second order central difference scheme in the space direction for continuous systems, and NSFD scheme is also applied, which leads to the proposal of the discrete system with time delay.
\item[$\bullet$]
Based on the Lyapunov function dependent on historical information, the global properties of the proposed discrete systems are investigated, including the existence and uniqueness of global positive solutions, the existence and global stability of equilibrium points, which are consistent with the corresponding continuous systems.
\item[$\bullet$]
The Lyapunov function constructed in this paper can guarantee the $\alpha$-robustness of L1 scheme, and the stability of two equilibrium points are the same.
\item[$\bullet$]
The influence of the theoretical results on disease transmission are verified by numerical simulation.
\end{itemize}

The organization of this paper is as follows. Some preliminaries are introduced in Section $2$, including definition and basic properties of L1 discrete scheme. In Section $3$, by using the second order central difference scheme and L1 scheme, the discrete system of the fractional-order reaction-diffusion SIR epidemic model is derived, and then the positivity, boundedness and global stability of two equilibrium points are established. Numerical simulations are presented to illustrate theoretical results in Section $4$. Finally, a brief discussion is given in Section $5$.

\section{Preliminaries }

This section begins with some definitions and results.

\textbf{Definition 2.1. \cite{3}} A Gamma function of $\alpha$ ($\mathfrak{R}(\alpha)>0$) is defined by:
\[\Gamma(\alpha)= \int_{0}^\infty{x^{\alpha-1}e^{-x}dx}.\]

\textbf{Definition 2.2. \cite{3}} Riemann-Liouville (R-L) fractional-order integral of order $\alpha $ ($n-1< \alpha< n)$ for a function $f(t)\in \mathbb{R}$ is defined by
\[I_{0+}^\alpha f\left( t \right)= \frac{1}{{\Gamma
\left( { \alpha } \right)}}\int_{t_0}^t \frac{{f}(s)}{{{{\left( {t - s}
\right)}^{1-\alpha} }}}ds,\]
where $\Gamma\left( \cdot \right)$ is the Gamma function.

\textbf{Definition 2.3. \cite{3}} Caputo fractional-order derivative of order $\alpha $ ($n-1< \alpha< n)$ for a function $f(t)\in \mathbb{C}^n[0,+\infty]$ is defined by
\[^{C}_{0}D^{\alpha}_{t} f\left( t \right)= \frac{1}{{\Gamma
\left( { n-\alpha } \right)}}\int_{0}^t \frac{{f^{(n)}}(s)}{{{{\left( {t - s}
\right)}^{\alpha-n+1} }}}ds,\]
where $\Gamma\left( \cdot \right)$ is the Gamma function and $f^{(n)}(s)=\frac{d^n f(s)}{ds^n}$. In particular, Caputo fractional-order derivative with order $0<\alpha<1$ is given by:
\[^{C}_{0}D^{\alpha}_{t} f\left( t \right)= \frac{1}{{\Gamma
\left( { 1-\alpha } \right)}}\int_{0}^t \frac{f'(s)}{{{{\left( {t - s}
\right)}^{\alpha} }}}ds.\]

\textbf{Definition 2.4. \cite{2006A}} The L1 discrete scheme of Caputo fractional derivatives for order $\alpha$ ($0<\alpha<1$) is defined by:
\begin{flalign}
\begin{aligned}
^{C}_{0}D^{\alpha}_{t_n}f(t_n)&=\frac{1}{\Gamma(1-\alpha)}\sum_{j=0}^{n-1}\int_{t_j}^{t_{j+1}}\frac{f'(\eta)}{(t_n-\eta)^\alpha}d\eta\\
&\approx \frac{(\Delta t)^{-\alpha}}{\Gamma(2-\alpha)}\sum_{j=0}^{n-1}b_j^{(1-\alpha)}[f(t_{n-j})-f(t_{n-j-1})]\\
&=\frac{(\Delta t)^{-\alpha}}{\Gamma(2-\alpha)}[f(t_n)-b_{n-1}^{(1-\alpha)}f(t_0)-\sum_{j=1}^{n-1}(b_{j-1}^{(1-\alpha)}-b_j^{(1-\alpha)})f(t_{n-j})]\\
&\triangleq \delta_n^\alpha f(t_n),\nonumber
\end{aligned}
\end{flalign}
where $b_j^{(1-\alpha)}=(j+1)^{1-\alpha}-j^{1-\alpha}$ and $\Delta t=t_{j+1}-t_j$ for $j=1,2,....$.

\textbf{Remark 2.1.} If the fractional order $\alpha=1$, the coefficients $b_0^{(1-\alpha)}=1$ and $b_j^{(1-\alpha)}=0$ for $j\ge 1$. Then for $n\ge 1$, the following equation holds:
\begin{flalign}
\begin{aligned}
\lim_{\alpha\to 1^-}\delta_n^\alpha f(t_n)=\frac{(\Delta t)^{-\alpha}}{\Gamma(1)}[f(t_n)-(b_{0}^{(1-\alpha)}-b_1^{(1-\alpha)})f(t_{n-1})]=\frac{1}{\tau}(f(t_n)-f(t_{n-1})).\nonumber
\end{aligned}
\end{flalign}
Then the above L1 discrete scheme is standard forward Euler difference scheme as follows:
\[\Delta f(t_n)=\frac{f(t_{n})-f(t_{n-1})}{\Delta t}.\]


\textbf{Definition 2.5. \cite{2017Stability} (Asymptotic stability of discrete systems)}  Consider the following nonautonomous system:
\[\delta_n^{\alpha}x_n=f(x_n).\]
Let $x_n=0$ be an equilibrium point of the above system. \\
(1) The equilibrium point $x_n=0$ is said to be stable if for all $\epsilon>0$, there exists a $\delta=\delta(t_0,\epsilon)>0$ such that if
\[||x_0||<\delta,~then~||x_n||<\epsilon,\]
for $n\ge 0$.\\
(2) The equilibrium point $x_n=0$ is asymptotically if there exists a $\delta=\delta(t_0)>0$, $||x_0||<\delta$ implies $\lim_{n\to\infty}x_n=0$.

\textbf{Remark 2.2} It can be found from Definition 2.6 that the asymptotic stability of equilibrium point for discrete system is consistent with that of continuous system \cite{2010Stability}.

\textbf{Definition 2.6. \cite{2002Nonlinear}} Let $\varphi:R^+\to R^+$ be a continuous function. $\varphi(r)$ is said to belong to class-K if it is strictly increasing and $\varphi(0) = 0$.

\textbf{Definition 2.7 \cite{2001ode}} Let $V(x):U\to R^n$ ($U\subset R^n$) be continuously differentiable function and $V(x^*)=0$. If all $x\in U$ satisfies $V(x)>0$ except $x=x^*$, $V(x)$ is said to be positive definite.


\textbf{Lemma 2.1. \cite{2020Barrier}} If $0<\alpha<1$, $b_k^{(\alpha)}=(k+1)^{\alpha}-k^{\alpha}$ satisfies the following qualities:\\
(1) $b_k^{(\alpha)}>b_{k+1}^{(\alpha)}$, $k=0,1,2,\cdots$;\\
(2) $b_0^{(\alpha)}=1$, $b_k^{(\alpha)}>0$, $k=0,1,2,\cdots$;\\
(3) there exists $C_1>0$ and $C_2>0$ such that $C_1k^{1-\alpha}\le (b_k^{(\alpha)})^{-1}\le C_2 k^{1-\alpha}$. 

In the fractional calculus, Cruz et al. proposed a useful inequality for Lyapunov function \cite{2015Volterra}. We extend this result to the following discrete case.

\textbf{Lemma 2.2.} Consider the positive sequence $x_n$, then for all $n\ge 1$, the following inequality holds:
\[\delta_n^\alpha(x_n-1-lnx_n)\le (1-\frac{1}{x_n})\delta_n^\alpha x_n.\]

\textbf{Proof:} Let define the difference scheme $\Delta_j (lnx_{n-j})=lnx_{n-j}-lnx_{n-j-1}$ and $\Delta_j (x_{n-j})=x_{n-j}-x_{n-j-1}$ for $j=1,2,...,n$. According to $\Delta_j (lnx_{n})=0$ and $\Delta_j (1)=0$, the following equation holds:
\begin{flalign}
\begin{aligned}
\delta_n^\alpha(lnx_n)-\frac{1}{x_n}\delta_n^\alpha x_n&=\frac{(\Delta t)^{-\alpha}}{\Gamma(2-\alpha)}\sum_{j=0}^{n-1}b_j^{(1-\alpha)}[\Delta_j (lnx_{n-j})-\frac{1}{x_n}\Delta_j (x_{n-j})]\\
&=\frac{(\Delta t)^{-\alpha}}{\Gamma(2-\alpha)}\sum_{j=0}^{n-1}b_j^{(1-\alpha)}\Delta_j (lnx_{n-j}-lnx_{n}-\frac{1}{x_n}x_{n-j}+1)\\
&=\frac{(\Delta t)^{-\alpha}}{\Gamma(2-\alpha)}[b_0^{(1-\alpha)}(lnx_{n}-lnx_{n}-1+1)-b_{n-1}^{(1-\alpha)}(lnx_{0}-lnx_{n}-\frac{x_0}{x_n}+1)\\
&~~~~-\sum_{j=0}^{n-1}(b_{j-1}^{(1-\alpha)}-b_{j}^{(1-\alpha)})(lnx_{n-j}-lnx_{n}+\frac{x_{n-j}}{x_n}+1)]\\
&=\frac{(\Delta t)^{-\alpha}}{\Gamma(2-\alpha)}[-b_{n-1}^{(1-\alpha)}(ln\frac{x_{0}}{x_{n}}+1-\frac{x_0}{x_n})-\sum_{j=0}^{n-1}(b_{j-1}^{(1-\alpha)}-b_{j}^{(1-\alpha)})(ln\frac{x_{n-j}}{x_{n}}+1-\frac{x_{n-j}}{x_n})].\nonumber
\end{aligned}
\end{flalign}
Obviously, $ln\frac{x_{i}}{x_{j}}+1-\frac{x_i}{x_j}\le 0$ for all $i,j=1,2,...$, so Eq.(1) holds as follows: \[\delta_n^\alpha(lnx_n)-\frac{1}{x_n}\delta_n^\alpha x_n\ge 0.\]
Then it can be yields the following conclusion:
\begin{flalign}
\begin{aligned}
\delta_n^\alpha(x_n-1-lnx_n)-(1-\frac{1}{x_n})\delta_n^\alpha x_n=-\delta_n^\alpha lnx_n +\frac{1}{x_n}\delta_n^\alpha x_n\le 0,\nonumber
\end{aligned}
\end{flalign}
that is,
\[\delta_n^\alpha(x_n-1-lnx_n)\le (1-\frac{1}{x_n})\delta_n^\alpha x_n.\]
The proof of Lemma 2.2 is completed. \hfill$\square$


\textbf{Lemma 2.3. \cite{2001ode}} Let $V(x):U\to R^n$ ($U\subset R^n$) be continuously differentiable function. Then $V(x)$ be a positive definite function if and only if there exist class-K functions $\varphi_1(x)$ and $\varphi_2(x)$ that satisfy the following equation:
\[\varphi_1(x)\le V(x)\le \varphi_2(x).\]

Lyapunovs direct method is an effective tool to analyze the stability of nonlinear systems without solving its state equations. In the following, we present the following theorem of a discrete fractional Lyapunov direct method with L1 scheme.

\textbf{Theorem 2.1 (Discrete fractional Lyapunov method)} Consider the following autonomous systems:
\begin{flalign}
\begin{aligned}\delta_n^{\alpha}x_n=l(x_n),
\end{aligned}
\end{flalign}
where $\alpha\in(0,1)$, $\delta_n^{\alpha}$ is the discrete scheme of Caputo fractional derivatives,  and $l(x_n):R\to R$ be a function such that there exists a unique solution of system (1) and $l(0)=0$. If there exists a positive define function $V(x_n)$ satisfying the following autonomous inequality:
\[\delta_n^{\alpha}V(x_n)\le -W(x_n),\]
where $W(x_n)$ is a positive define function and $W(x_n)=0$ if and only if $x_n=0$. Then $x_n=0$ is globally asymptotically stable.

\textbf{Proof} It is obvious from $\delta_n^{\alpha}V(x_n)\le 0$ that one has $V(x_n)\le V(x_0)$, which implies $x_n=0$ is stable. According to the positive define functions $W(x_n)$ and $V(x_n)$ and Lemma 2.3, there exists class-K functions $\varphi_1(x_n)$, $\varphi_2(x_n)$, $\tilde{\varphi}_1(x_n)$ and $\tilde{\varphi}_2(x_n)$ satisfied
\[{\varphi}_{1}(x_n)\le W(x_n)\le {\varphi}_{2}(x_n),\]
and
\[\tilde{{\varphi}}_{1}(x_n)\le V(x_n)\le \tilde{{\varphi}}_{2}(x_n),\]
then it is obvious that
\[\delta_n^{\alpha}V(x_n)\le -\tilde{W}(V(x_n)),\]
where $\tilde{W}(V(x_n))={\varphi}_1(\tilde{\varphi}_2^{-1}(V(x_n))$ is a monotonically increasing function. The construction auxiliary function satisfies $\delta_n^{\alpha}\tilde{V}(x_n)= -\tilde{W}(\tilde{V}(x_n))$, which implies $V(x_n)\le \tilde{V}(x_n)$ ($n=1,2,...$). Therefore it is easy to see from \cite{2019Complete} that $\tilde{V}(x_n)$ is a monotonically decreasing function. So one has that there exists a positive constant $M$ satisfying $\lim_{n\to \infty} \tilde{V}(x_n)=M$, which implies $\lim_{n\to \infty} \delta_n^{\alpha}\tilde{V}(x_n)=0$. Thus it can be deduced that $\lim_{n\to \infty} \tilde{W}(\tilde{V}(x_n))=0$. Furthermore, based on the the monotonicity of $\tilde{W}(V(x_n))$, it is straightforward to show $\lim_{n\to \infty} x_n=0$.  Then $x_n=0$ is globally asymptotically stable. \hfill$\square$

\textbf{Remark 2.3. (Discrete Lyapunov method)} According to Remarks 2.1 and 2.2, when $\alpha=1$, one has $\delta_n^\alpha x_n=\frac{x_{n}-x_{n-1}}{\tau}$. So it is easy to see that Theorem 2.1 is valid for $\alpha=1$, that is:\\
Consider the following autonomous systems:
\begin{flalign}
\begin{aligned}
x_{n}=x_{n-1}+\tau l(x_n),\nonumber
\end{aligned}
\end{flalign}
and $l(0)=0$. If there exists a positive define function $V(x_n)$ satisfying the following autonomous inequality:
\[\Delta V(x_n)\le -W(x_n),\]
where $W(x_n)$ is a positive define function and $W(x_n)=0$ if and only if $x_n=0$. Then $x_n=0$ is globally asymptotically stable.

\section{System Dscription}
Over the past few decades, an intensive effort has been put into developing theoretical models for anomalous diffusion that cannot be modeled as standard Brownian motion \cite{2000From,2000The}. Henry et al. indicated that the signature of the above anomalous diffusion is that the mean square displacement of the diffusing species $<(\Delta x)^2>$ scales as a nonlinear power law in time. i.e. $<(\Delta x)^2>\sim t^\alpha$. Meanwhile, faster than linear scaling $(\alpha>1)$ is referred to as superdiffusion and slower than linear scaling $(0<\alpha<1)$ is referred to as subdiffusion \cite{2006Anomalous}. Furthermore, Meerschaert et al. proposed that a power law waiting time distribution leads to fractional-order derivative with the same order \cite{2012Stochastic}. There has a number of paper considering fractional-order anomalous diffusion models of infectious diseases. Particular, Zhen et al. investigated a class of fractional-order SIR epidemic model with reaction-diffusion by the general incidence rate in \cite{2021luzhenzhen} as follows:
\begin{flalign}
\begin{aligned}
\left\{ \begin{array}{l}
_0^CD_t^\alpha {S} = {d_{1}}\Delta {S} + \lambda- \gamma{S} -{\beta}S{f}(I)\\
_0^CD_t^\alpha {I} = {d_{2}}\Delta {I} + \beta Sf(I)-(\mu + r){I}\\
_0^CD_t^\alpha {R} = {d_{3}}\Delta {R} + {r}{I} - \delta{R}\\
\end{array} \right.
\end{aligned}
\end{flalign}
where S, I and R represent the number of susceptible, infected and recovered individuals respectively; $_0^CD_t^\alpha$ implies Caputo fractional-order operator ($0<\alpha< 1$); $\Delta {\rm{ = }}\frac{{{\partial ^2}}}{{\partial {x^2}}}$ $(x\in \Omega \subset \mathbb{R}^2)$ denotes the Laplace operator; $\gamma$ and $\delta$ imply the nature death rates of $S$ and $R$, respectively; $\mu$ denotes the disease-related death rates of $I$; $\lambda$ represents the recruitment rate of the total population; $r$ implies the recovery rate of the infected individuals; $d_{i}$ $(k=1,2,3)$ denotes the diffusion rate of $S$, $I$ and $R$, respectively; $\beta$ is the infection rate of $S$ infected by $I$. Besides, in the process of infectious disease transmission, susceptible individuals contact infected individuals and become infected with a certain probability. There are lot of evidences that the incidence rate is an important tool to describe this process \cite{anderson1992infectious,bailey1975mathematical}. But in reality, detailed information about infectious diseases is often difficult to obtain, as such information can change with the surrounding environment and gradually decrease with an individual's perception. Therefore, the general incidence rate $f(I)$ will be considered. Meanwhile, the generalized incidence rate $f(I)$ satisfies the following assumptions in this paper:
\begin{flalign}
\begin{aligned}
&\textbf{(A1)}:~f(0)=0~and~f(I)>0~for~I>0;\\
&\textbf{(A2)}:~f'(I)>0~and~(\frac{f}{I})'\le 0~for~I\ge 0;\\
&\textbf{(A3)}:~f(I)\le If'(0). \nonumber
\end{aligned}
\end{flalign}
Biologically, assumptions ${\bf (A1)}-{\bf(A3)}$ mean that the disease transmission rates are monotonically increasing, but subject to saturation effects, such as $f(I)=I$ and $f(I)=\frac{I}{w+I}$, where the saturation effect $w>0$ \cite{2020Global}. In this paper, suppose there is no external input in system (2), that is no population movement on the boundary $\partial \Omega$, which implies the Neumann boundary is considered:
\begin{flalign}
\begin{aligned}
\frac{\partial S}{\partial \upsilon}= \frac{\partial I}{\partial \upsilon } = \frac{\partial R}{\partial \upsilon } = 0,x \in \partial \Omega,~t \ge 0,
\end{aligned}
\end{flalign}
and the initial condition is as follows:
\begin{flalign}
\begin{aligned}
({S}(0,x),{I}(0,x),{R}(0,x))= ({\phi_{1}}(x),{\phi_{2}}(x),{\phi_{3}(x)}),~x \in \Omega.
\end{aligned}
\end{flalign}
For sake of convenience, we denote system (2)-(4) by system (2) throughout this paper.

Furthermore, the following lemma is proposed without proof:

\textbf{Lemma 3.1. \cite{2021luzhenzhen}} There exists a unique nonnegative solution $(S(x,t),I(x,t),R(x,t))$ of system (2), which is also ultimately bounded for any given initial function $\phi(x)=(\phi_{1}(x),\phi_{2}(x),\phi_{3}(x))\neq(0,0,0)$.

\textbf{Lemma 3.2. \cite{2021luzhenzhen}} Under the reproduction number $R_0={(\frac{{{\beta}S^0{f'(0)}}}{{{\mu} + {r}}})}\le 1$ and assumptions ${\bf (A1)}-{\bf(A3)}$ holding, there exists the unique disease-free equilibrium point $E_0=(\frac{\lambda}{\gamma},0,0)$ of system (2), which is globally asymptotically stable.

\textbf{Lemma 3.3. \cite{2021luzhenzhen}} Under $R_0>1$ and assumptions ${\bf (A1)}-{\bf(A3)}$ holding, system (2) is uniform persistence, that is for any initial value $\phi(x)=(\phi_{1}(x),\phi_{2}(x),\phi_{3}(x))\neq(0,0,0)$, the solution $(S(t,\phi),I(t,\phi),R(t,\phi))$ satisfies
\[\mathop {\liminf}\limits_{t \to \infty } {S}(t,\phi ) \ge w ,~\mathop {\liminf}\limits_{t \to \infty } {I}(t,\phi ) \ge w,~\mathop {\liminf}\limits_{t \to \infty } {R}(t,\phi ) \ge w\]
where $w$ is a positive constant. Furthermore, system (2) has at least one endemic equilibrium $E^*=(S^*,I^*,R^*)$ satisfying
\[\left\{ \begin{array}{l}
{\lambda} -{\gamma}{S^*} - {\beta }{S^*}{f}({I^*})=0 \\
 {{\beta}{S^*}{f}({I})}  - ({\mu } + {r}){I^*}=0\\
 rI^*-\delta R^*=0,
\end{array} \right.\]
which is globally asymptotically stable.

In the following section, we dedicate ourselves to the study of the discrete counterpart of the continuous system (2).

\subsection{Preliminary results}
Consider system (2) in the spital domain $\Omega=[a,b]\times [a,b]$. To discretize system (2), we let $\Delta x=(b-a)/M$ be the space stepsize generating $M$ equal subintervals over the domain and $\Delta t$ be the time stepsize. Then at each mesh point $(x_n,t_k)$ with $x_n=a+n\Delta t$ and $t_k=k\Delta t$, we denote approximations of $S(x_n,t_k)$, $I(x_n,t_k)$ and $R(x_n,t_k)$ by $S_n^k$, $I_n^k$ and $R_n^k$, respectively. By applying L1 nonstandard finite difference scheme and second order center difference scheme, the discretization system of system (2) is given by:
\begin{flalign}
\begin{aligned}
\frac{(\Delta t)^{-\alpha}}{\Gamma(2-\alpha)}[S_n^{k+1}-b_k^{(1-\alpha)}S_n^0-\sum_{j=1}^k(b_{j-1}^{(1-\alpha)}-b_{j}^{(1-\alpha)})S_n^{k+1-j}]&=d_1\frac{(S_{n-1}^{k+1}-2S_n^{k+1}+S_{n+1}^{k+1})}{(\Delta x)^2}+\lambda-\beta S_n^{k+1}f(I_n^k)-\gamma S_n^{k+1},\\
\frac{(\Delta t)^{-\alpha}}{\Gamma(2-\alpha)}[I_n^{k+1}-b_k^{(1-\alpha)}I_n^0-\sum_{j=1}^k(b_{j-1}^{(1-\alpha)}-b_{j}^{(1-\alpha)})I_n^{k+1-j}]&=d_2\frac{(I_{n-1}^{k+1}-2I_n^{k+1}+I_{n+1}^{k+1})}{(\Delta x)^2}+\beta S_n^{k+1}f(I_n^k)-(\mu+r) I_n^{k+1},\\
\frac{(\Delta t)^{-\alpha}}{\Gamma(2-\alpha)}[R_n^{k+1}-b_k^{(1-\alpha)}R_n^0-\sum_{j=1}^k(b_{j-1}^{(1-\alpha)}-b_{j}^{(1-\alpha)})R_n^{k+1-j}]&=d_3\frac{(R_{n-1}^{k+1}-2R_n^{k+1}+R_{n+1}^{k+1})}{(\Delta x)^2}+rI_n^{k+1}-\delta R_n^{k+1}.
\end{aligned}
\end{flalign}
Rearrangement system (5) yields the following discrete system:

\begin{flalign}
\begin{aligned}
S_n^{k+1}-r_1(S_{n-1}^{k+1}-2S_n^{k+1}+S_{n+1}^{k+1})&=b_k^{(1-\alpha)}S_n^0+\sum_{j=1}^k(b_{j-1}^{(1-\alpha)}-b_{j}^{(1-\alpha)})S_n^{k+1-j}+g(\lambda-\beta S_n^{k+1}f(I_n^k)-\gamma S_n^{k+1}),\\
I_n^{k+1}-r_2(I_{n-1}^{k+1}-2I_n^{k+1}+I_{n+1}^{k+1})&=b_k^{(1-\alpha)}I_n^0+\sum_{j=1}^k(b_{j-1}^{(1-\alpha)}-b_{j}^{(1-\alpha)})I_n^{k+1-j}+ g(\beta S_n^{k+1}f(I_n^k)-(\mu+r) I_n^{k+1}),\\
R_n^{k+1}-r_3(R_{n-1}^{k+1}-2R_n^{k+1}+R_{n+1}^{k+1})&=b_k^{(1-\alpha)}R_n^0+\sum_{j=1}^k(b_{j-1}^{(1-\alpha)}-b_{j}^{(1-\alpha)})R_n^{k+1-j}+g(rI_n^{k+1}-\delta R_n^{k+1}),
\end{aligned}
\end{flalign}
where $n\in \{0,1,2,\cdots,M\}$, $k\in N$, $r_i=\frac{d_i\Gamma(2-\alpha) (\Delta t)^\alpha}{(\Delta x)^2}$ $(i=1,2,3)$ and $g=\Gamma(2-\alpha)(\Delta t)^\alpha$, with the initial condition:
\begin{flalign}
\begin{aligned}
S_n^0=\phi_1(x_n)>0,~I_n^0=\phi_2(x_n)>0,~R_n^0=\phi_3(x_n)>0.
\end{aligned}
\end{flalign}
Without losing generality, system (6) will be studied next. Meanwhile, first order center difference scheme is applied to the Neumann boundary condition as follows:
\begin{flalign}
\begin{aligned}
S_{-1}^k=S_0^k,~I_{-1}^k=I_0^k,~R_{-1}^k=R_0^k.
\end{aligned}
\end{flalign}
Similarly, we denote system (6)-(8) by system (6) throughout this paper.

\textbf{Remark 3.1.} Particularly, according to Remark 2.1, when $\alpha=1$, system (6) can be further simplified as follows:
\begin{flalign}
\begin{aligned}
\frac{S_n^{k+1}-S_n^{k}}{\Delta t}&=d_1\frac{(S_{n-1}^{k+1}-2S_n^{k+1}+S_{n+1}^{k+1})}{(\Delta x)^2}+\lambda-\beta S_n^{k+1}f(I_n^k)-\gamma S_n^{k+1},\\
\frac{I_n^{k+1}-I_n^{k}}{\Delta t}&=d_2\frac{(I_{n-1}^{k+1}-2I_n^{k+1}+I_{n+1}^{k+1})}{(\Delta x)^2}+\beta S_n^{k+1}f(I_n^k)-(\mu+r) I_n^{k+1},\\
\frac{R_n^{k+1}-R_n^{k}}{\Delta t}&=d_3\frac{(R_{n-1}^{k+1}-2R_n^{k+1}+R_{n+1}^{k+1})}{(\Delta x)^2}+rI_n^{k+1}-\delta R_n^{k+1},
\end{aligned}
\end{flalign}
which is consistent with the integer order discretization model in \cite{2018Global}. Therefore, the fractional-order
continuous system discretized by L1 scheme and second order central difference scheme is consistent with the integer-order continuous system discretized by forward difference scheme when the fractional-order $\alpha=1$, which implies L1 scheme is $\alpha$-robust for the discretization model (6) when $\alpha\to 1^-$.

\textbf{Remark 3.2.} It is worth noting that system (6) considers a discrete time delay system, but unlike ordinary delay systems \cite{2015Lyapunov}, there is no delay at the initial value, and the time delay of system (6) starts from $t_1$, so there is no need to guarantee the delay of the initial conditions (7).

In the following, the non-negative bounded solution of system (6) will be investigated.

\textbf{Theorem 3.1.} For any $\Delta t>0$ and $\Delta x>0$, the solution of system (6) is non-negative and bounded for any $k\in N$.

\textbf{Proof:} For system (6), it can be written in matrix form as follows:
\begin{flalign}
\begin{aligned}
\left\{ \begin{array}{l}
A^k S^{k+1}=b_k^{(1-\alpha)}S^0+\sum_{j=1}^k (b_{j-1}^{(1-\alpha)}-b_j^{(1-\alpha)})S^{k+1-j}+g\lambda,\\
B^k I^{k+1}=b_k^{(1-\alpha)}I^0+\sum_{j=1}^k (b_{j-1}^{(1-\alpha)}-b_j^{(1-\alpha)})I^{k+1-j}+g\beta S^{k+1}f(I^k),\\
C^k R^{k+1}=b_k^{(1-\alpha)}R^0+\sum_{j=1}^k (b_{j-1}^{(1-\alpha)}-b_j^{(1-\alpha)})R^{k+1-j}+grI^{k+1},\\
\end{array} \right.
\end{aligned}
\end{flalign}
where $S^k=(S_0^k,S_1^k,\cdots,S_M^k)$, $I^k=(I_0^k,I_1^k,\cdots,I_M^k)$, $R^k=(R_0^k,R_1^k,\cdots,R_M^k)$ ($k=0,1,2,\dots,N$)
\[A^k=\left(
                 \begin{array}{ccccccc}
                   a_0^k & a & 0 & \cdots & 0 & 0 & 0 \\
                   a & a_1^k & a & \cdots & 0 & 0 & 0\\
                   \vdots &  \vdots &  \vdots &  \vdots &  \vdots &  \vdots &  \vdots \\
                   0 &  0&  0 &\cdots &a &a_{M-1}^k &a\\
                   0 & 0 & 0 & \cdots & 0& a & a_M^k \\
                 \end{array}
               \right),
~B^k=\left(
                 \begin{array}{ccccccc}
                   b_0^k & b & 0 & \cdots & 0 & 0 & 0 \\
                   b & b_1^k & b & \cdots & 0 & 0 & 0\\
                   \vdots &  \vdots &  \vdots &  \vdots &  \vdots &  \vdots &  \vdots \\
                   0 &  0&  0 &\cdots &b &b_{M-1}^k &b\\
                   0 & 0 & 0 & \cdots &0 & b & b_M^k \\
                 \end{array}
               \right),\]
\[C^k=\left(
                 \begin{array}{ccccccc}
                   c_0^k & c & 0 & \cdots & 0 & 0 & 0 \\
                   c & c_1^k & c & \cdots & 0 & 0 & 0\\
                   \vdots &  \vdots &  \vdots &  \vdots &  \vdots &  \vdots &  \vdots \\
                   0 &  0&  0 &\cdots &c &c_{M-1}^k &c\\
                   0 & 0 & 0 & \cdots &0 & c & c_M^k \\
                 \end{array}
               \right),\]
and $a_0^k=1+r_1+g\beta f(I_0^k)+g\gamma$, $a_i^k=1+2r_1+g\beta f(I_i^k)+g\gamma$ $(i=1,2,\cdots,M-1)$, $a_M^k=1+r_1+g\beta f(I_M^k)+g\gamma$, $a=-r_1$, $b_0^k=1+r_2+g(\mu+r)$, $b_i^k=1+2r_2+g(\mu+r)$ $(i=1,2,\cdots,M-1)$, $b_M^k=1+r_2+g(\mu+r)$, $b=-r_2$, $c_0^k=1+r_3+g\delta$, $c_i^k=1+2r_3+g\delta$ $(i=1,2,\cdots,M-1)$, $c_M^k=1+r_3+g\delta$ and $c=-r_3$. It is clear that $A^k$, $B^k$ and $C^k$ are a strictly diagonally matrix, which imply $(A^k)^{-1}$, $(B^k)^{-1}$ and $(C^k)^{-1}$ exists and are positive definite matrices. Thus, system (10) is equivalent to:
\begin{flalign}
\begin{aligned}
\left\{ \begin{array}{l}
S^{k+1}=(A^k)^{-1} [b_k^{(1-\alpha)}S^0+\sum_{j=1}^k (b_{j-1}^{(1-\alpha)}-b_j^{(1-\alpha)})S^{k+1-j}+g\lambda],\\
I^{k+1}=(B^k)^{-1} [b_k^{(1-\alpha)}I^0+\sum_{j=1}^k (b_{j-1}^{(1-\alpha)}-b_j^{(1-\alpha)})I^{k+1-j}+g\beta S^{k+1}f(I^k)],\\
R^{k+1}=(C^k)^{-1} [b_k^{(1-\alpha)}R^0+\sum_{j=1}^k (b_{j-1}^{(1-\alpha)}-b_j^{(1-\alpha)})R^{k+1-j}+grI^{k+1}].\\
\end{array} \right.
\end{aligned}
\end{flalign}
implying that there exists a unique solution of system (6). Since all parameters in system (11) are positive with $b_{j-1}^{(1-\alpha)}-b_j^{(1-\alpha)}>0$, it is easy to see that the solution remain non-negative for all $k\in N$.

Next the boundedness of the solution will be investigated. Define a sequence $\{G^k\}$ as follows:
\[G^k=\sum_{n=0}^M(S_n^k+I_n^k+R_n^k).\]
It follows from system (6) that one has
\[G^{k+1}\le b_k^{(1-\alpha)}G^0+\sum_{j=1}^k (b_{j-1}^{(1-\alpha)}-b_j^{(1-\alpha)})G^{k+1-j}+g\lambda-g\tilde{\mu}G^{k+1}.\]
where $\tilde{\mu}=\min\{\mu,\gamma,\delta\}$. Particular, when $k=0$, it can be arrived at
\[G^1\le \frac{G^0+g\lambda}{1+g\tilde{\mu}}.\]
Similarly, when $k=1$, one has
\[G^2\le \frac{G^0+g\lambda+G^1}{1+g\tilde{\mu}}\le \frac{G^0+g\lambda}{1+g\tilde{\mu}}+\frac{G^0+g\lambda}{(1+g\tilde{\mu})^2}\le(G^0+g\lambda)\sum_{j=1}^2\frac{j}{(1+g\tilde{\mu})^j}.\]
and when $k=2$, it is easy to see
\[G^3\le \frac{G^0+g\lambda+G^1+G^2}{1+g\tilde{\mu}}\le \frac{G^0+g\lambda}{1+g\tilde{\mu}}+\frac{2(G^0+g\lambda)}{(1+g\tilde{\mu})^2}+\frac{(G^0+g\lambda)}{(1+g\tilde{\mu})^3}\le(G^0+g\lambda)\sum_{j=1}^3\frac{j}{(1+g\tilde{\mu})^j}.\] So it can be obtained by mathematical induction that when $k=n$, one has
\[G^{n+1}\le (G^0+g\lambda)\sum_{j=1}^{n+1}\frac{j}{(1+g\tilde{\mu})^j}.\]
As we know $\sum_{n=1}^\infty nx^n=\frac{x}{(1-x)^2}$ ($-1<x<1$), so it can be yielded from the above equation that
\[\lim_{n\to\infty}G^{n+1}\le(G^0+g\lambda)\sum_{j=1}^{\infty}\frac{j}{(1+g\tilde{\mu})^j}=\frac{(g\lambda+G^0)(1+g\tilde{\mu})}{(g\tilde{\mu})^2}.\]
So $S_n^k$, $I_n^k$ and $R_n^k$ are uniformly bounded eventually.
\hfill$\square$

\subsection{Equilibrium points}
This section is devoted to the investigation of the existence of the disease-free equilibrium point and the endemic equilibrium point.

\textbf{Theorem 3.2.} System (6) has a unique disease-free equilibrium point $E_0=(\frac{\lambda}{\gamma},0,0)$, which is consistent with the disease-free equilibrium point of the corresponding continuous system (2). And when the reproduction number $R_0=\beta\frac{\frac{\lambda}{\gamma}f'(0)}{\mu+r}>1$ and assumption {\textbf{(A1)-(A3)}} hold, system (6) exists a unique endemic equilibrium point $E^*=(S^*,I^*,R^*)$ satisfying:
\begin{flalign}
\begin{aligned}
\left\{ \begin{array}{l}
\lambda-\beta S^*f(I^*)-\gamma S^*=0,\\
\beta S^*f(I^*)-(\mu+r)I^*=0,\\
rI^*-\delta R^*=0,\\
\end{array} \right.\nonumber
\end{aligned}
\end{flalign}
which is also consistent with the corresponding continuous system (2).

\textbf{Proof:} A disease-free equilibrium point is a steady state solution of system (6) with all infected variables equal to zero, that is $I_n^k=R_n^k=0$ $(n=1,2,...M;~k\in N)$.  So the disease-free equilibrium point $E_0$ of system (6) is satisfied the followed equation:
\begin{flalign}
\begin{aligned}
S_{0}&=b_k^{(1-\alpha)}S_0+\sum_{j=1}^k (b_{j-1}^{(1-\alpha)}-b_j^{(1-\alpha)})S_0+g\lambda-g\gamma S_0,\\
&=b_k^{(1-\alpha)}S_0+(b_0^{(1-\alpha)}-b_k^{(1-\alpha)})S_0+g\lambda-g\gamma S_0,\nonumber
\end{aligned}
\end{flalign}
So $g\lambda-g\gamma S_0=0$, that is $S_0=\frac{\lambda}{\gamma}$, which is consistent with the corresponding continuous system (2). Similarly, if exists, the endemic equilibrium point $E^*$ satisfies the following equation:
\begin{flalign}
\begin{aligned}
\left\{ \begin{array}{l}
\lambda-\beta S^*f(I^*)-\gamma S^*=0,\\
\beta S^*f(I^*)-(\mu+r)I^*=0,\\
rI^*-\delta R^*=0.\\
\end{array} \right.\nonumber
\end{aligned}
\end{flalign}
Let $g(I)=\beta \frac{\lambda-(\mu+r)I}{\gamma}f(I)-(\mu+r)I$, then $g(0)=0$, $g(\frac{\lambda}{\mu+r})=-\lambda<0$ and
\[g'(0)=\beta\frac{\lambda f'(0)}{\gamma}-(\mu+r)>0,\]
as the reproduction number $R_0=\beta\frac{\frac{\lambda}{\gamma}f'(0)}{\mu+r}> 1$, that is $g(I)$ decreases monotonically at $I=0$, then due to the continuity of $g(I)$, there exists a constant $I^*\ in (0,\frac{\lambda}{\mu+r})$ make $g(I^*)=0$.  Therefore, when the reproduction number $R_0>1$, there exists the endemic equilibrium point $E^*=(S^*,I^*,R^*)$ of system (6). The uniqueness of $E^*$ will be analyzed in the following. Using the fact that
\begin{flalign}
\begin{aligned}
g'(I^*)&=\beta\frac{\lambda-(\mu+r)I^*}{\gamma}f'(I^*)-\beta\frac{\mu+r}{\gamma}f(I^*)-(\mu+r)\\
&=\beta S^*f'(I^*)-\beta S^*f(I^*)-\beta S^*\frac{f(I^*)}{I^*}\\
&=\beta\frac{S^*}{I^*}[I^*f'(I^*)-f(I^*)]-S^*f(I^*).
\end{aligned}
\end{flalign}
According to assumption {\bf(A3)}, we have $If'(I)\le f(I)$, which implies that $g'(I^*)<0$. If there exists the second positive equilibrium $E^{**}=(S^{**},I^{**},R^{**})$, one has $g'(I^{**})>0$. But this is a contradiction with Eq.(12). Therefore, there exists a unique equilibrium point of system (6). \hfill$\square$

\subsection{Global stability of equilibrium points}
In this subsection, the global stability of two equilibrium points will be investigated. The following lemma will be considered first before the main results are given:

\textbf{Lemma 3.4.} Considering the sequence $u^k$ and $h^k$ satisfy $u^{k+1}=\sum_{i=1}^{k+1}w_{k+2-i}h^i$, \[u^0=\frac{(w_1b_k^{(1-\alpha)}+w_2b_{k-1}^{(1-\alpha)}+...+w_{k+1}b_0^{(1-\alpha)})h^0}{b_k^{(1-\alpha)}},\]
and $w_k$ are positive constants. Then it can be given as follows:
\[\delta_n^\alpha u^{k+1}=\sum_{i=1}^{k+1}w_{k+2-i}\delta_n^\alpha h^i.\]

\textbf{Proof:} Calculating L1 scheme of $u^{k+1}$ ($i=1,2,...,k+1$) satisfied:
\begin{flalign}
\begin{aligned}
\delta_n^\alpha u^{k+1}=\frac{\tau^{-\alpha}}{\Gamma(2-\alpha)}[u^{k+1}-b_{k}^{(1-\alpha)}u^0-\sum_{j=1}^{k}(b_{j-1}^{(1-\alpha)}-b_{j}^{(1-\alpha)})u^{k+1-j}].
\end{aligned}
\end{flalign}
Substituting the function $u^i$ ($i=0,1,...,k+1$) into Eq.(13) can be yielded:
\begin{flalign}
\begin{aligned}
\delta_n^\alpha u^{k+1}&=\frac{\tau^{-\alpha}}{\Gamma(2-\alpha)}[\sum_{i=1}^{k+1}w_{k+2-i}h^i-\sum_{i=1}^{k+1}w_{k+2-i}b_{i-1}^{(1-\alpha)}h^0-\sum_{j=1}^{k}(b_{j-1}^{(1-\alpha)}-b_{j}^{(1-\alpha)})\sum_{i=1}^{k+1-j}w_{k+2-i-j}h^i]\\
&=\frac{\tau^{-\alpha}}{\Gamma(2-\alpha)}[\sum_{i=1}^{k+1}w_{k+2-i}h^i-\sum_{i=1}^{k+1}w_{k+2-i}b_{i-1}^{(1-\alpha)}h^0-\sum_{j=1}^{k}(b_{j-1}^{(1-\alpha)}-b_{j}^{(1-\alpha)})\sum_{i=j}^{k+1}w_{k+2-i}h^{i-j}]\\
&=\frac{\tau^{-\alpha}}{\Gamma(2-\alpha)}[\sum_{i=1}^{k+1}w_{k+2-i}h^i-\sum_{i=1}^{k+1}b_{i-1}^{(1-\alpha)}w_{k+2-i}h^0-\sum_{i=1}^{k+1}w_{k+2-i}\sum_{j=1}^{k}(b_{j-1}^{(1-\alpha)}-b_{j}^{(1-\alpha)})h^{i-j}]\\
&=\frac{\tau^{-\alpha}}{\Gamma(2-\alpha)}\sum_{i=1}^{k+1}w_{k+2-i}[h^i-b_{i-1}^{(1-\alpha)}h^0-\sum_{j=1}^{k}(b_{j-1}^{(1-\alpha)}-b_{j}^{(1-\alpha)})h^{i-j}]\\
&=\sum_{i=1}^{k+1}w_{k+2-i}\delta_n^\alpha h^{i}.\nonumber
\end{aligned}
\end{flalign}
Therefore, it is easy to see $\delta_n^\alpha u^{k+1}=\sum_{i=1}^{k+1}w_{k+2-i}\delta_n^\alpha h^i$. \hfill$\square$

\textbf{Theorem 3.3.} Supposed that assumptions $\textbf{(A1)-(A3)}$ hold. For any $\Delta t>0$ and $\Delta x>0$, if the reproduction number $R_0=\beta\frac{\lambda f'(0)}{\gamma (\mu+r)}\le 1$, the disease-free equilibrium point $E_0$ of system (6) is globally asymptotically stable.
\textbf{Proof:} Define a discrete Lyapunov function as follows:
\[W^{k+1}=\sum_{i=1}^{k+1}w_{k+2-i}V^{i}+\beta S_0f'(0)\sum_{n=0}^MI_n^{k+1},\]
where
\[V^{i}=\sum_{n=0}^M(S_n^{i}-S_0-S_0ln\frac{S_n^{i}}{S_0})+\sum_{n=0}^MI_n^{i},\]
for $i=1,...,k+1$, $V^0=\sum_{n=0}^M(S_n^{0}-S_0-S_0ln\frac{S_n^{0}}{S_0})+\sum_{n=0}^MI_n^{0}$, $W^0=\frac{(w_1b_k^{(1-\alpha)}+w_2b_{k-1}^{(1-\alpha)}+...+w_{k+1}b_0^{(1-\alpha)})V^0}{b_k^{(1-\alpha)}}$ and the parameter $w_i$ are positive constants to be determined later. Then $W^{k+1}$ is a positive definite function and $W^{k+1}=0$ if and only if $S_n^{i}=S_0$ and $I^{i}_n=0$ ($i=1,2,...,k+1$). Calculating L1 scheme of $V^{i}$ ($i=1,2,...,k+1$) satisfied:
\begin{flalign}
\begin{aligned}
\delta_n^\alpha V^{i}&\le \sum_{n=0}^M (1-\frac{S_0}{S_n^{i}})\delta_n^\alpha S_n^{i}+\sum_{n=0}^M \delta_n^\alpha I_n^{i}\\
&=\sum_{n=0}^M [(1-\frac{S_0}{S_n^{i}})(\lambda-\beta S_n^{i}f(I_n^{i-1})-\gamma S_n^{i})+(\beta S_n^{i}f(I_n^{i-1})-(\mu+r)I_n^{i})]+\Delta L_1^{i}+\Delta L_2^{i}\\
&=\sum_{n=0}^M [-\frac{\gamma(S_n^{i}-S_0)^2}{S_n^{i}}-\beta S_n^{i}f(I_n^{i-1})+\beta S_0f(I_n^{i-1})+\beta S_n^{i}f(I_n^{i-1})-(\mu+r)I_n^{i}]+\Delta L_1^{i}+\Delta L_2^{i}\\
&\le \sum_{n=0}^M[ -\frac{\gamma(S_n^{i}-S_0)^2}{S_n^{i}}+\beta S_0f'(0)I_n^{i-1}-(\mu+r)I_n^{i}]+\Delta L_1^{i}+\Delta L_2^{i}\\
&=\sum_{n=0}^M [-\frac{\gamma(S_n^{i}-S_0)^2}{S_n^{i}}+(\beta S_0f'(0)-(\mu+r))I_n^{i}+\beta S_0f'(0)(I_n^{i-1}-I_n^{i})]+\Delta L_1^{i}+\Delta L_2^{i},
\end{aligned}
\end{flalign}
where
\begin{flalign}
\begin{aligned}
\Delta L_1^{i}+\Delta L_2^{i}=\sum_{n=0}^M \frac{1}{(\Delta x)^2}[d_1(1-\frac{S_0}{S_n^{i}})(S_{n+1}^{i}-2S_n^{i}+S_{n-1}^{i})+d_2(I_{n+1}^{i}-2I_n^{i}+I_{n-1}^{i})].
\end{aligned}
\end{flalign}
Since it is important to note that
Eq.(15) in this case, can be written as
\begin{flalign}
\begin{aligned}
&~~~~~\Delta L_1^{i}+\Delta L_2^{i}\\
&=\sum_{n=0}^M \frac{1}{(\Delta x)^2}[d_1(S_{n+1}^{i}-2S_n^{i}+S_{n-1}^{i})-S_0d_1(\frac{S_{n+1}^{i}}{S_n^{i}}-2+\frac{S_{n}^{i}}{S_{n+1}^{i}})+d_2(I_{n+1}^{i}-2I_n^{i}+I_{n-1}^{i})]\\
&\le \frac{1}{(\Delta x)^2}[d_1(S_{M+1}^{i}-2S_M^{k+1}+S_{-1}^{i}-S_0^{i})-S_0d_1(\frac{S_{-1}^{i}}{S_0^{i}}-2+\frac{S_{M+1}^{i}}{S_{M}^{i}})+d_2(I_{M+1}^{i}-2I_M^{i}+I_{-1}^{i}-I_0^{k+1})]\\
&=0.
\end{aligned}
\end{flalign}
Therefore, substitute Eq.(16) into Eq.(14) is given by:
\begin{flalign}
\begin{aligned}
\delta_n^\alpha V^{i}\le \sum_{n=0}^M [-\frac{\gamma(S_n^{i}-S_0)^2}{S_n^{i}}+(\beta S_0f'(0)-(\mu+r))I_n^{i}+\beta S_0f'(0)(I_n^{i-1}-I_n^{i})].
\end{aligned}
\end{flalign}
for $i=1,2,...,k+1$. Also taking $\delta_n^\alpha I_n^k$ to obtain the following equation:
\begin{flalign}
\begin{aligned}
\delta_n^\alpha I_n^{k+1}=\frac{\tau^{-\alpha}}{\Gamma(2-\alpha)}[I_n^{k+1}-b_{k}^{(1-\alpha)}I_n^0-\sum_{j=1}^{k}(b_{j-1}^{(1-\alpha)}-b_{j}^{(1-\alpha)})I_n^{k+1-j}].
\end{aligned}
\end{flalign}
Let $\Delta_n^i=\beta S_0f'(0)(I_n^{i-1}-I_n^{i})$ for $i=1,2,...,k+1$. Then calculating the sum of $\Delta_n^i$ multiplied by $w_{k+2-i}$ and adding Eq. (18) can be yielded:
\begin{flalign}
\begin{aligned}
&~~~~\sum_{n=0}^M [\sum_{i=1}^{k+1}w_{k+2-i}\Delta_n^i+\beta S_0f'(0)\delta_n^\alpha I_n^{k+1}]\\
&=\sum_{n=0}^M\sum_{i=1}^{k+1}w_{k+2-i}\beta S_0f'(0)(I_n^{i-1}-I_n^{i})+\beta S_0f'(0)\sum_{n=0}^M\delta_n^\alpha I_n^{k+1}\\
&=\beta S_0f'(0)\{\sum_{n=0}^M[(-w_1+\frac{\tau^{-\alpha}}{\Gamma(2-\alpha)})I^{k+1}_n+(w_1-w_2-\frac{\tau^{-\alpha}}{\Gamma(2-\alpha)}(b_{0}^{(1-\alpha)}-b_{1}^{(1-\alpha)}))I^k_n\\
&~~~~+(w_2-w_3-\frac{\tau^{-\alpha}}{\Gamma(2-\alpha)}(b_{1}^{(1-\alpha)}-b_{2}^{(1-\alpha)}))I^{k-1}_n+...\\
&~~~~+(w_k-w_{k+1}-\frac{\tau^{-\alpha}}{\Gamma(2-\alpha)}(b_{k-1}^{(1-\alpha)}-b_{k}^{(1-\alpha)}))I_n^1+(w_{k+1}-\frac{\tau^{-\alpha}}{\Gamma(2-\alpha)}b_{k}^{(1-\alpha)})I^0]\}.
\end{aligned}
\end{flalign}
Take $w_i= \frac{{(\Delta t)}^{-\alpha}}{\Gamma(2-\alpha)}b_{i-1}^{(1-\alpha)}$ ($i\ge 1$), then it can be yielded
\[\sum_{n=0}^M(\sum_{i=1}^{k+1}w_{k+2-i}\Delta_n^i+\beta S_0f'(0)\delta_n^\alpha I_n^k)=0.\]
According to Lemma 3.4 and $w_i\le \frac{{(\Delta t)}^{-\alpha}}{\Gamma(2-\alpha)}$, calculating the fractional-order difference of $W^{k+1}$ satisfies:
\begin{flalign}
\begin{aligned}
\delta_n^\alpha W^{k+1}&=\sum_{i=1}^{k+1}w_{k+2-i}\delta_n^\alpha V^{i}+\beta S_0f'(0)\sum_{n=0}^M\delta_n^\alpha I_n^{k+1}\\
&\le \sum_{i=1}^{k+1}w_{k+2-i} \sum_{n=0}^M (1-\frac{S_0}{S_n^{i}})\delta_n^\alpha S_n^{i}+\sum_{i=1}^{k+1}w_{k+2-i} \sum_{n=0}^M \delta_n^\alpha I_n^{i}+\sum_{n=0}^M\delta_n^\alpha I_n^{k+1}\\
&\le\sum_{i=1}^{k+1}w_{k+2-i} \sum_{n=0}^M [-\frac{\gamma(S_n^{i}-S_0)^2}{S_n^{i}}+(\beta S_0f'(0)-(\mu+r))I_n^{i}]+[\sum_{i=1}^{k+1}w_{k+2-i}\Delta_n^i+\beta S_0f'(0)\sum_{n=0}^M \delta_n^\alpha I_n^{k+1}]\\
&\le \sum_{i=1}^{k+1}w_{k+2-i} \sum_{n=0}^M [-\frac{\gamma(S_n^{i}-S_0)^2}{S_n^{i}}+(\beta S_0f'(0)-(\mu+r))I_n^{i}].
\end{aligned}
\end{flalign}
Define $\tilde{V}^{k+1}$ satisfying the following equation:
\[\tilde{V}^{k+1}=\frac{{(\Delta t)}^{-\alpha}}{\Gamma(2-\alpha)}\sum_{i=1}^{k+1}w_{k+2-i}\sum_{n=0}^M [\frac{\gamma(S_n^{i}-S_0)^2}{S_n^{i}}-(\mu+r)(R_0-1)I_n^{i}].\]
then $\tilde{V}^{k+1}$ is a positive define function if the reproduction number $R_0\le 1$ and $\tilde{V}^{k+1} = 0$ if and only if $S^{k+1}_n = S_0$ and $I_n^{k+1}=0$. Then according to Theorem 2.1, one has
\[\lim_{i\to \infty}S_n^i=S_0, ~\lim_{i\to\infty}I_n^i=0,~\lim_{i\to\infty}R_n^i=0, ~for~n=1,2,...,M.\]
Then, if the reproduction number $R_0=\beta\frac{\lambda f'(0)}{\gamma (\mu+r)}\le 1$, the disease-free equilibrium point $E_0=(S_0,0,0)$ of system (6) is globally asymptotically stable.  \hfill$\square$

Next, we discuss the global stability of the endemic equilibrium $E^*$ when the reproduction number $R_0>1$.

\textbf{Theorem 3.4.} Supposed that assumptions $\textbf{(A1)-(A3)}$ hold. For any $\Delta t>0$ and $\Delta x>0$, if the reproduction number $R_0>1$, the endemic equilibrium point $E^*=(S^*,I^*,R^*)$ of system (6) is globally asymptotically stable.

\textbf{Proof:} Define the discrete Lyapunov function as follows:
\[W^{k+1}=\sum_{i=1}^{k+1}w_{k+2-i}(H_1^{i}+H_2^i)+\beta S^*f(I^*)H_3^{k+1},\]
where
\begin{flalign}
\begin{aligned}
H_1^{i}&=\sum_{n=0}^M (S_n^{i}-S^*-S^*ln\frac{S_n^{i}}{S^*}),\\
H_2^{i}&=\sum_{n=0}^M (I_n^{i}-I^*-I^*ln\frac{I_n^{i}}{I^*}),~i=1,2,...,k+1,\nonumber
\end{aligned}
\end{flalign}
and
\[H_3^{k+1}=\sum_{n=0}^M (\frac{I_n^{k+1}}{I^*}-1-ln\frac{I_n^{k+1}}{I^*})\triangleq\sum_{n=0}^M \Phi(\frac{I_n^{k+1}}{I^*}),\]
where $\Phi(w)=w-1-lnw$. In addition, let $W^0=\frac{(w_1b_k^{(1-\alpha)}+w_2b_{k-1}^{(1-\alpha)}+...+w_{k+1}b_0^{(1-\alpha)})(H_1^0+H_2^0)}{b_k^{(1-\alpha)}}$ where $H_1^{0}=\sum_{n=0}^M (S_n^{0}-S^*-S^*ln\frac{S_n^{0}}{S^*})$ and $H_2^{0}=\sum_{n=0}^M (I_n^{0}-I^*-I^*ln\frac{I_n^{0}}{I^*})$. The parameter $w_i>0$ are positive constants to be determined later. Obviously, $H_1^{i}\ge 0$, $H_2^{i}\ge 0$ ($i=1,2,...,k+1$) and $H_3^{k}\ge 0$ with the equality holds if and only if $S_n^k=S^*$ and $I_n^k=I^*$ for all $n\in\{0,1,2,...,M\}$ and $k\in N$, respectively. Calculating L1 scheme of $H_3^{k+1}$ can be yielded:
\begin{flalign}
\begin{aligned}
\delta_n^\alpha H_3^{k+1}=\sum_{n=0}^M \frac{\tau^{-\alpha}}{\Gamma(2-\alpha)}[\Phi(\frac{I_n^{k+1}}{I^*})-b_{k}^{(1-\alpha)} \Phi(\frac{I_n^{0}}{I^*})-\sum_{j=1}^{k}(b_{j-1}^{(1-\alpha)}-b_{j}^{(1-\alpha)}) \Phi(\frac{I_n^{k-j}}{I^*})].
\end{aligned}
\end{flalign}
Furthermore, it can be arrived at:
\begin{flalign}
\begin{aligned}
\delta_n^\alpha H_1^{i}&\le \sum_{n=0}^M(1-\frac{S^*}{S_n^{i}})\delta_n^\alpha S_n^{i}\\
&=\sum_{n=0}^M(1-\frac{S^*}{S_n^{i}})[\lambda-\beta S_n^{i}f(I_n^{i-1})-\gamma S_n^{i}]+\Delta H_1^{i}\\
&=\sum_{n=0}^M(1-\frac{S^*}{S_n^{i}})[\beta S^*f(I^*)+\gamma S^*-\beta S_n^{i}f(I_n^{i-1})-\gamma S_n^{i}]+\Delta H_1^{i}\\
&=\sum_{n=0}^M \frac{-\gamma(S_n^{i}-S^*)^2}{S_n^{i}}+\sum_{n=0}^M (1-\frac{S^*}{S_n^{i}})[\beta S^*f(I^*)-\beta S_n^{i}f(I_n^{i-1})]+\Delta H_1^{i},
\end{aligned}
\end{flalign}
and
\begin{flalign}
\begin{aligned}
\delta_n^\alpha H_2^{i}&\le \sum_{n=0}^M(1-\frac{I^*}{I_n^{i}})\delta_n^\alpha I_n^{i}\\
&=\sum_{n=0}^M(1-\frac{I^*}{I_n^{i}})[\beta S_n^{i}f(I_n^{i-1})-(\mu+r)I_n^{i}]+\Delta H_2^{i},
\end{aligned}
\end{flalign}
where
\[\Delta H_1^{i}=\sum_{n=0}^M \frac{1}{(\Delta x)^2}d_1(1-\frac{S^*}{S_n^{i}})(S_{n+1}^{i}-2S_n^{i}+S_{n-1}^{i}),\]
and
\[\Delta H_2^{i}=\sum_{n=0}^M \frac{1}{(\Delta x)^2}d_2(1-\frac{I^*}{I_n^{i}})(I_{n+1}^{i}-2I_n^{i}+I_{n-1}^{i}).\]
Similar to Theorem 3.3 and applying
\[\frac{S_{n+1}^{i}}{S_n^{i}}-2+\frac{S_{n}^{i}}{S_{n+1}^{i}}\ge 0,~\frac{I_{n+1}^{i}}{I_n^{k+1}}-2+\frac{I_{n}^{i}}{I_{n+1}^{i}}\ge 0,\]
it is straightforward to show that
\[\Delta H_1^{i}\le 0, ~\Delta H_2^{i}\le 0,\]
for $i=1,2,...,k+1$. Therefore, define $H^{i}=H_1^{i}+H_2^{i}$ and it is obvious from Eq.(22) and Eq.(23) that
\begin{flalign}
\begin{aligned}
\delta_n^\alpha H^{i}&\le \sum_{n=0}^M \frac{-\gamma(S_n^{i}-S^*)^2}{S_n^{i}}+\sum_{n=0}^M (1-\frac{S^*}{S_n^{i}})[\beta S^*f(I^*)-\beta S_n^{i}f(I_n^{i-1})]\\
&~~~~+\sum_{n=0}^M(1-\frac{I^*}{I_n^{i}})[\beta S_n^{i}f(I_n^{i-1})-(\mu+r)I_n^{i}]\\
&=\sum_{n=0}^M \frac{-\gamma(S_n^{i}-S^*)^2}{S_n^{i}}+\sum_{n=0}^M \beta S^*f(I^*)[2-\frac{S^*}{S_n^{i}}-\frac{I^*}{I_n^{i}}+\frac{f(I_n^{i-1})}{f(I^*)}-\frac{f(I_n^{i-1})S_n^{i}I^*}{f(I^*)S^*I_n^{i}}].
\end{aligned}
\end{flalign}
By some calculations, one has
\begin{flalign}
\begin{aligned}
&~~~~~~~2-\frac{S^*}{S_n^{i}}-\frac{I^*}{I_n^{i}}+\frac{f(I_n^{i-1})}{f(I^*)}-\frac{f(I_n^{i-1})S_n^{i}I^*}{f(I^*)S^*I_n^{i}}\\
&=(1-\frac{S^*}{S_n^{i}}+ln\frac{S^*}{S_n^{i}})+(1-\frac{f(I_n^{i-1})S_n^{i}I^*}{f(I^*)S^*I_n^{i}}+ln\frac{f(I_n^{i-1})S_n^{i}I^*}{f(I^*)S^*I_n^{i}})\\
&~~~~~~+(-\frac{I_n^{i}}{I^*}+ln\frac{I_n^{i}}{I^*}+\frac{I_n^{i-1}}{I^*}-ln\frac{I_n^{i-1}}{I^*})+(1-\frac{I_n^{i-1} f(I^*)}{I^* f(I_n^{i-1})}+ln\frac{I_n^{i-1} f(I^*)}{I^* f(I_n^{i-1})})-1+\frac{I_n^{i-1}f(I^*)}{I^*f(I_n^{i-1})}+\frac{f(I_n^{i-1})}{f(I^*)}-\frac{I_n^{i-1}}{I^*}\\
&=-\Phi(\frac{S^*}{S_n^{i}})-\Phi(\frac{f(I_n^{i-1})S_n^{i}I^*}{f(I^*)S^*I_n^{i}})-\Phi(\frac{I_n^{i-1}f(I^*)}{I^*f(I_n^{i-1})})+\Phi(\frac{I_n^{i-1}}{I^*})-\Phi(\frac{I_n^{i}}{I^*})+\frac{I_n^*}{f(I^*)f(I_n^{i-1})}[f(I_n^{i-1})-f(I^*)][\frac{f(I_n^{i-1})}{I_n^{i-1}}-\frac{f(I^*)}{I^*}].
\end{aligned}
\end{flalign}
Based on assumptions \textbf{(A2)} and \textbf{(A3)}, Eq.(25) further simplify to give:
\begin{flalign}
\begin{aligned}
2-\frac{S^*}{S_n^{i}}-\frac{I^*}{I_n^{i}}+\frac{f(I_n^{i-1})}{f(I^*)}-\frac{f(I_n^{i-1})S_n^{i}I^*}{f(I^*)S^*I_n^{i}}
\le-\Phi(\frac{S^*}{S_n^{i}})-\Phi(\frac{f(I_n^{i-1})S_n^{i}I^*}{f(I^*)S^*I_n^{i}})-\Phi(\frac{I_n^{i-1}f(I^*)}{I^*f(I_n^{i-1})})+\Phi(\frac{I_n^{i-1}}{I^*})-\Phi(\frac{I_n^{i}}{I^*}).
\end{aligned}
\end{flalign}
Then this can be substituted into Eq.(24), giving:
\begin{flalign}
\begin{aligned}
\delta_n^\alpha H^{i}&\le \sum_{n=0}^M \frac{-\gamma(S_n^{i}-S^*)^2}{S_n^{i}}+\sum_{n=0}^M \beta S^*f(I^*)[-\Phi(\frac{S^*}{S_n^{i}})-\Phi(\frac{f(I_n^{i-1})S_n^{i}I^*}{f(I^*)S^*I_n^{i}})-\Phi(\frac{I_n^{i-1}f(I^*)}{I^*f(I_n^{i-1})})+\Phi(\frac{I_n^{i-1}}{I^*})-\Phi(\frac{I_n^{i}}{I^*})].
\end{aligned}
\end{flalign}
Let $\Delta_i=\sum_n^M [\Phi(\frac{I_n^{i-1}}{I^*})-\Phi(\frac{I_n^{i}}{I^*})]$ and calculate the following equation satisfying:
\begin{flalign}
\begin{aligned}
&~~~~\sum_{n=0}^M [\sum_{i=1}^{k+1}w_{k+2-i}\beta S^*f(I^*)\Delta_n^i+\beta S^*f(I^*)\delta_n^\alpha \Phi(\frac{I_n^{k+1}}{I^*})]\\
&=\beta S^*f(I^*)[\sum_{n=0}^M\sum_{i=1}^{k+1}w_{k+2-i}\beta S_0f'(0)(I_n^{i-1}-I_n^{i})+\sum_{n=0}^M\delta_n^\alpha \Phi(\frac{I_n^{k+1}}{I^*})]\\
&=\beta S^*f(I^*)\sum_{n=0}^M[(-w_1+\frac{\tau^{-\alpha}}{\Gamma(2-\alpha)})\Phi(\frac{I_n^{k+1}}{I^*})+(w_1-w_2-\frac{\tau^{-\alpha}}{\Gamma(2-\alpha)}(b_{0}^{(1-\alpha)}-b_{1}^{(1-\alpha)}))\Phi(\frac{I_n^{k}}{I^*})+...\\
&~~~~+(w_k-w_{k+1}-\frac{\tau^{-\alpha}}{\Gamma(2-\alpha)}(b_{k-1}^{(1-\alpha)}-b_{k}^{(1-\alpha)}))\Phi(\frac{I_n^{1}}{I^*})+(w_{k+1}-a\frac{\tau^{-\alpha}}{\Gamma(2-\alpha)}b_{k}^{(1-\alpha)})\Phi(\frac{I_n^{0}}{I^*})].
\end{aligned}
\end{flalign}
Take $w_i= \frac{{(\Delta t)}^{-\alpha}}{\Gamma(2-\alpha)}b_{i-1}^{(1-\alpha)}$ ($i\ge 1$), then the following equation holds:
\[\sum_{n=0}^M(\sum_{i=1}^{k+1}w_{k+2-i}\Delta_n^i+\beta S^*f(I^*)\delta_n^\alpha \Phi(\frac{I_n^{k+1}}{I^*}))=0.\]
Calculating the fractional-order difference of $W^{k+1}$ according to Lemma 3.4, satisfies
\[\delta_n^\alpha W^{k+1}=\sum_{i=1}^{k+1}w_{k+2-i}(\delta_n^\alpha H_1^{i}+\delta_n^\alpha H_2^i)+\beta S^*f(I^*)\delta_n^\alpha H_3^{k+1}.\]
Then substitute Eqs.(24), (27) and (28) yields the following equation:
\begin{flalign}
\begin{aligned}
\delta_n^\alpha W^{k+1}&\le \sum_{i=1}^{k+1}w_{k+2-i}\{\sum_{n=0}^M \frac{-\gamma(S_n^{i}-S^*)^2}{S_n^{i}}+\sum_{n=0}^M \beta S^*f(I^*)[-\Phi(\frac{S^*}{S_n^{i}})-\Phi(\frac{f(I_n^{i-1})S_n^{i}I^*}{f(I^*)S^*I_n^{i}})-\Phi(\frac{I_n^{i-1}f(I^*)}{I^*f(I_n^{i-1})})]\}.
\end{aligned}
\end{flalign}
Let define the following function:
\begin{flalign}
\begin{aligned}
\tilde{V}^{k+1}=\sum_{i=1}^{k+1}w_{k+2-i}\{\sum_{n=0}^M \frac{-\gamma(S_n^{i}-S^*)^2}{S_n^{i}}+\sum_{n=0}^M \beta S^*f(I^*)[-\Phi(\frac{S^*}{S_n^{i}})-\Phi(\frac{f(I_n^{i-1})S_n^{i}I^*}{f(I^*)S^*I_n^{i}})-\Phi(\frac{I_n^{i-1}f(I^*)}{I^*f(I_n^{i-1})})]\},\nonumber
\end{aligned}
\end{flalign}
then it is obvious $\delta_n^\alpha W^{k+1}\le -\tilde{V}^{k+1}$ and $\tilde{V}^{k+1}$ is a positive define function and $\tilde{V}^{k+1}=0$ if and only if $S_n^{k+1}=S^*$ and $I_n^{k+1}=I^*$. Thus, according to Theorem 2.1, it can be found that the following conclusion holds:
\[\lim_{k\to \infty}S_n^k=S^*, ~\lim_{k\to\infty}I_n^k=I^*,~\lim_{k\to\infty}R_n^k=R^*.\]
 Thus, the endemic equilibrium point $E^*=(S^*,I^*,R^*)$ is globally asymptotically stable. \hfill$\square$

\textbf{Remark 3.3.} It can be found that the results of Lemma 3.2, Lemma 3.3 and Theorem 3.3 and Theorem 3.4 are the same for the stability of two kinds of equilibrium points. In other words, L1 nonstandard finite difference scheme and second order central difference scheme can maintain the properties of the corresponding continuous system (2). The stability of equilibrium point of a numerical system represents that of corresponding continuous system.

\textbf{Remark 3.4.} Reviewing Eq.(19) again, when $\alpha=1$, one has $b_0^{(1-\alpha)}=1$ and $b_j^{(1-\alpha)}=0$ ($j\ge 1$), then the coefficient $w_k=0$ $(k>1)$, $w_1=\frac{1}{\Delta t}$. In this case, the following equation still holds:
\begin{flalign}
\begin{aligned}
&~~~~\sum_{n=0}^M (\sum_{i=1}^{k+1}w_{k+2-i}\Delta_n^i+\beta S_0f'(0)\delta_n^\alpha I_n^k)=0.\nonumber
\end{aligned}
\end{flalign}
Thus we can see that the conclusion in Theorem 3.3 is true for $\alpha=1$. In this case, the Lyapunov function is as follows:
\begin{flalign}
\begin{aligned}
W^{k+1}&=w_1[\sum_{n=0}^M(S_n^{k+1}-S_0-S_0ln\frac{S_n^{k+1}}{S_0})+\sum_{n=0}^MI_n^{k+1}]+\beta S_0f'(0)\sum_{n=0}^MI_n^{k+1}\\
&=\sum_{n=0}^M\frac{1}{\Delta t}[S_n^{k+1}-S_0-S_0ln\frac{S_n^{k+1}}{S_0}+(1+\beta S_0f'(0)\Delta t)I_n^{k+1}],\nonumber
\end{aligned}
\end{flalign}
which is consistence with \cite{2018Global}. Similarly, Eq.(28) is also true, then the conclusion is verified for Theorem 3.4. Meanwhile, $\alpha=1$ implies system (6) is simplified to system (9). Therefore, Theorem 3.3 and Theorem 3.4 hold for $\alpha=1$, which implies that L1 scheme is $\alpha-$robust in the investigation of Lyapunov functions of system (6), consistent with \cite{2019error}.

\textbf{Remark 3.5.} According to \cite{2018Global}, when the reproduction number $R_0\le 1$, the disease-free equilibrium point $E_0$ of system (9) is globally asymptotically stable, whereas the endemic equilibrium point $E^*$ of system (9) is globally asymptotically stable when $R_0>1$. This is consistent with system (6) for $\alpha=1$ about Theorem 3.3 and Theorem 3.4, which means L1 scheme is $\alpha$-robust in the stability of the system as $\alpha\to 1^-$.

\textbf{Remark 3.6.} In discussing of the disease-free equilibrium point, the Lyapunov function constructed in \cite{2018Global} is of the following form:
\[L^{k+1}=\sum_{n=0}^M\frac{1}{\Delta t}[S_n^{k+1}-S_0-S_0ln\frac{S_n^{k+1}}{S_0}+(1+\beta S_0f'(0)\Delta t)I_n^{k+1}],\]
which doesn't depend at all on the history information. However, the Lyapunov function constructed in this paper is strictly dependent on the information of before $k+1$:
\[W^{k+1}=\sum_{i=1}^{k+1}w_{k+2-i}(\sum_{n=0}^M(S_n^{i}-S_0-S_0ln\frac{S_n^{i}}{S_0})+\sum_{n=0}^MI_n^{i})+\beta S_0f'(0)\sum_{n=0}^MI_n^{k+1},\]
which is consistent with the property of discrete fractional-order operator, because L1 scheme of Caputo fractional-order operator is a memory operator. The same is true of the endemic equilibrium point.

\section{Numerical Simulation}
In this section, numerical simulations are given to illustrate the obtain results.

Consider the bilinear incidence rate $Sf(I)=SI$ as a example and $x\in R$, then system (2) becomes
\begin{flalign}
\begin{aligned}
\left\{ \begin{array}{l}
_{0}^CD^{\alpha}_{t}S(x,t)=d_1\frac{\partial^2 S}{\partial x^2}+\lambda-\beta S(x,t)I(x,t)-\gamma S,\\
_{0}^CD^{\alpha}_{t}I(x,t)=d_2\frac{\partial^2 I}{\partial x^2}+\beta S(x,t)I(x,t)-(\mu+r)I,\\
_{0}^CD^{\alpha}_{t}R(x,t)=d_3\frac{\partial^2 R}{\partial x^2}+rI-\delta R,\\
\end{array}\right.
\end{aligned}
\end{flalign}
Applying L1 nonstandard finite difference scheme and second order central difference scheme, the discretization of system (30) is given as follows:
\begin{flalign}
\begin{aligned}
S_n^{k+1}-r_1(S_{n-1}^{k+1}-2S_n^{k+1}+S_{n+1}^{k+1})&=b_k^{(1-\alpha)}S_n^0+\sum_{j=1}^k(b_{j-1}^{(1-\alpha)}-b_{j}^{(1-\alpha)})S_n^{k+1-j}+g(\lambda-\beta S_n^{k+1}I_n^k-\gamma S_n^{k+1}),\\
I_n^{k+1}-r_2(I_{n-1}^{k+1}-2I_n^{k+1}+I_{n+1}^{k+1})&=b_k^{(1-\alpha)}I_n^0+\sum_{j=1}^k(b_{j-1}^{(1-\alpha)}-b_{j}^{(1-\alpha)})I_n^{k+1-j}+ g(\beta S_n^{k+1}I_n^k-(\mu+r) I_n^{k+1}),\\
R_n^{k+1}-r_3(R_{n-1}^{k+1}-2R_n^{k+1}+R_{n+1}^{k+1})&=b_k^{(1-\alpha)}R_n^0+\sum_{j=1}^k(b_{j-1}^{(1-\alpha)}-b_{j}^{(1-\alpha)})R_n^{k+1-j}+g(rI_n^{k+1}-\delta R_n^{k+1}),\\
\end{aligned}
\end{flalign}
where $n\in \{0,1,2,\cdots,M\}$, $k\in N$, $r_i=\frac{d_i\Gamma(2-\alpha) (\Delta t)^\alpha}{(\Delta x)^2}$ $(i=1,2,3)$ and $g=\Gamma(2-\alpha)(\Delta t)^\alpha$.

In that follows, the initial value condition is taken as
\[S_n^0=0.5,~I_n^0=e^{-nh},~R_n^0=e^{-nh},\]
and the Neumann boundary condition is as follows:
\[S_{-1}^k=S_0^k,~I_{-1}^k=I_0^k,~R_{-1}^k=R_0^k.\]
The basic reproduction number of system (30) is $R_0=\frac{\beta \lambda}{\gamma(\mu+r)}$. System (31) possesses a disease-free equilibrium $E_0=(\frac{\lambda}{\mu},0,0)$. For $R_0>1$, system (31) also has an endemic equilibrium $E^*=(\frac{\mu+r}{\beta},\frac{\lambda}{\mu+r}-\frac{\gamma}{\beta},\frac{\lambda r}{(\mu+r)\delta}-\frac{r\gamma}{\beta\delta})$.

\begin{figure}
    \begin{center}
       \includegraphics[width=0.4\linewidth]{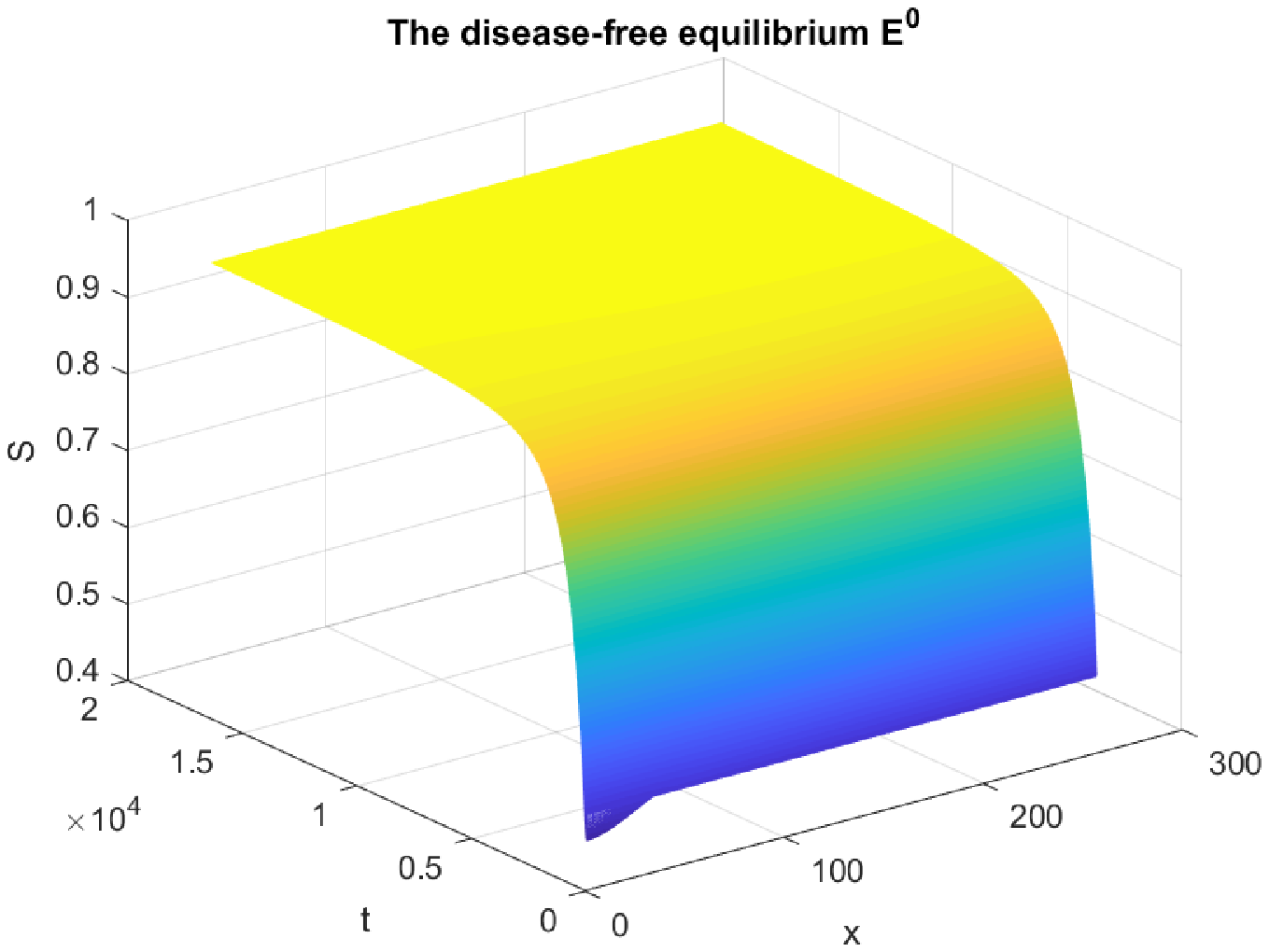}
       \includegraphics[width=0.4\linewidth]{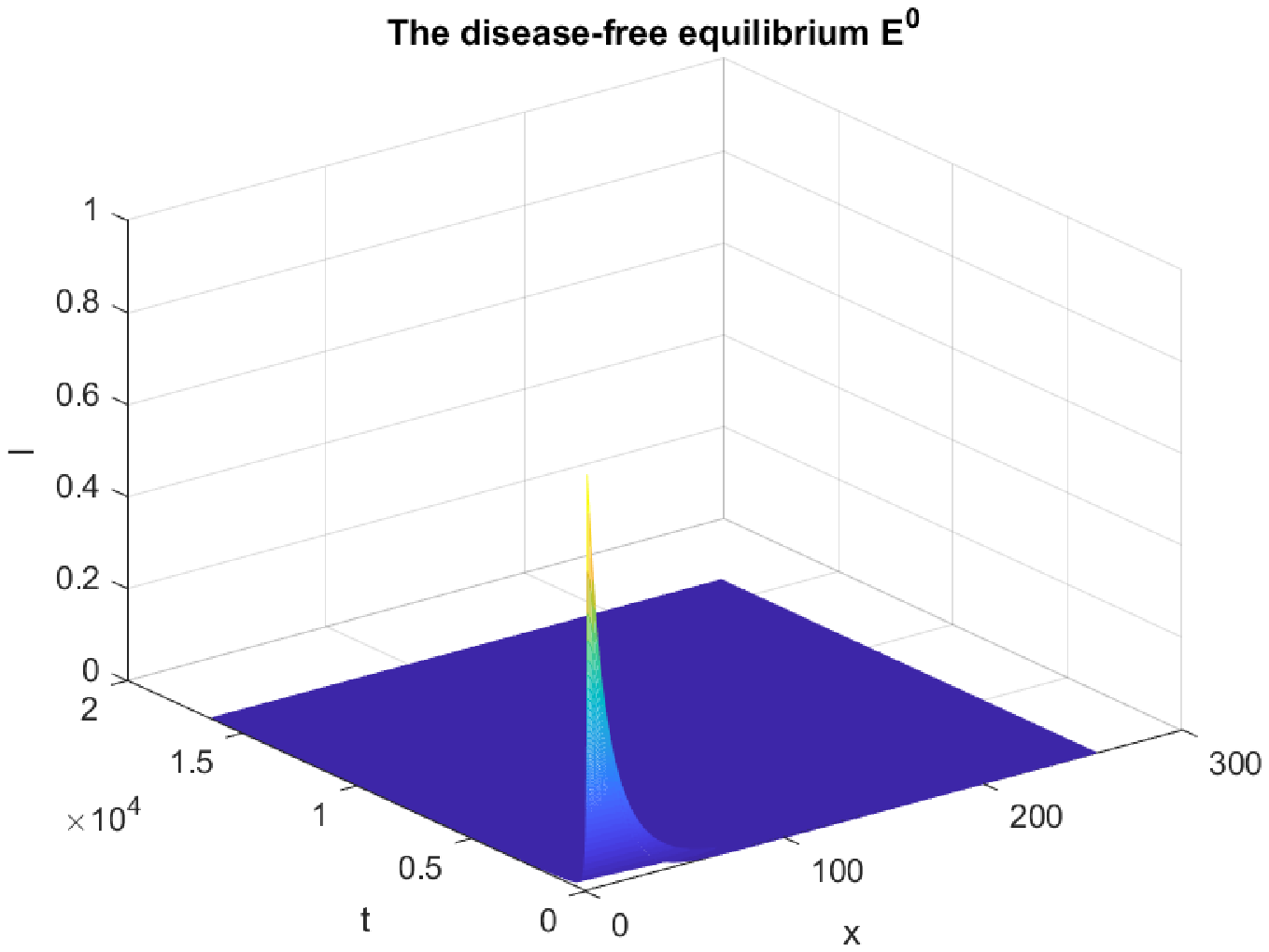}
       \includegraphics[width=0.4\linewidth]{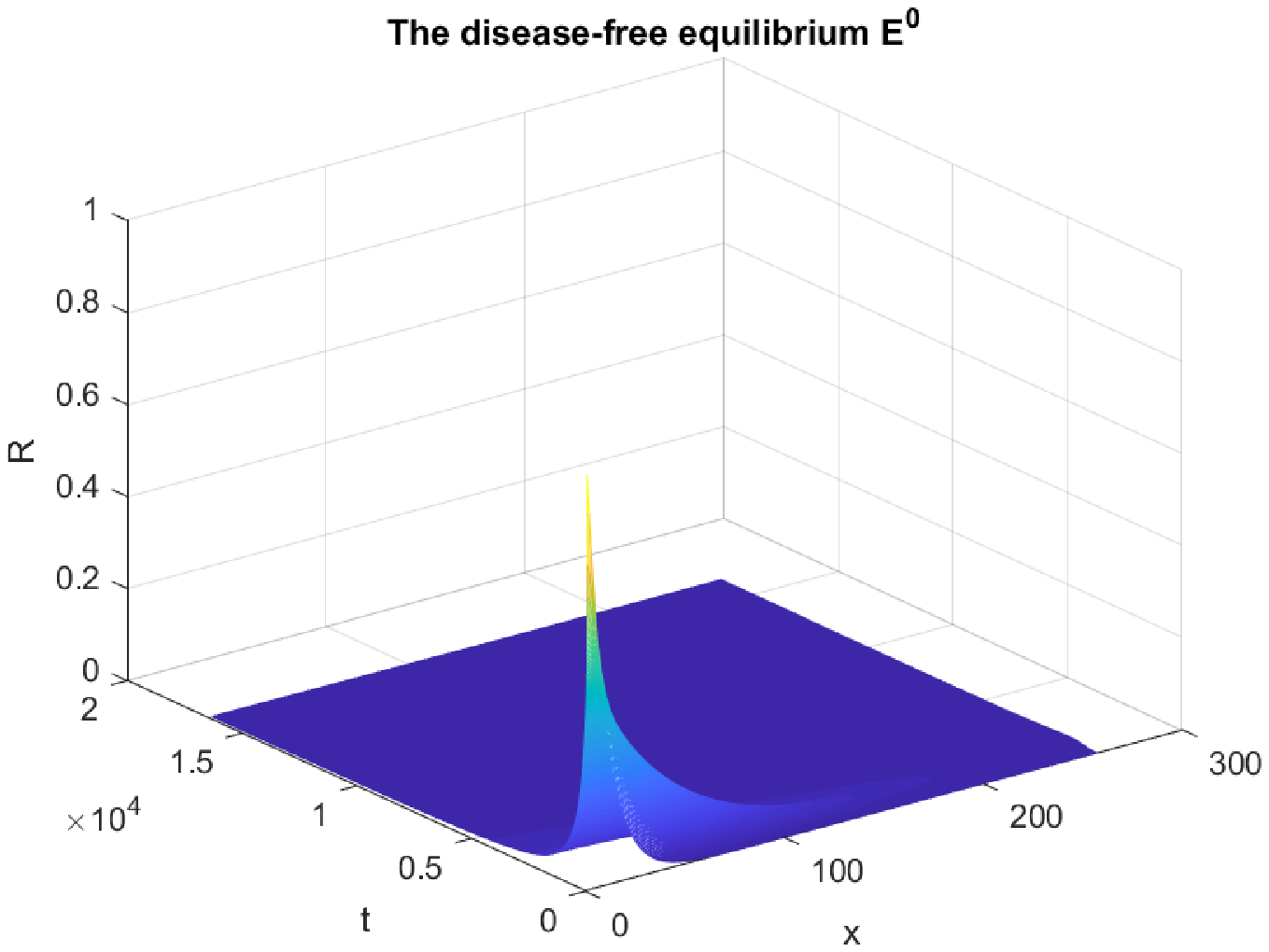}
      \caption{The disease free equilibrium $E_0$ is globally asymptotically stables.}
      \label{2}
    \end{center}
    \end{figure}

First let $\lambda=0.2$, $\beta=0.2144$, $\gamma=\delta=0.2$, $\mu=0.2$, $r=0.25$, $\alpha=0.8$, $d_1=d_2=d_3=1$, $\Delta t=0.1$ and $\Delta x=0.1$. In this case, the reproduction number $R_0=0.7238<1$ and the disease-free equilibrium point $E_0=(1,0,0)$. By Lemma 3.2 and Theorem 3.3, $E_0$ is globally asymptotically stable. This yield that the infectious disease is cleared, the infection dies out and the patient will be completely cured. Fig.~\ref{2} confirms these observations.

\begin{figure}
    \begin{center}
       \includegraphics[width=0.4\linewidth]{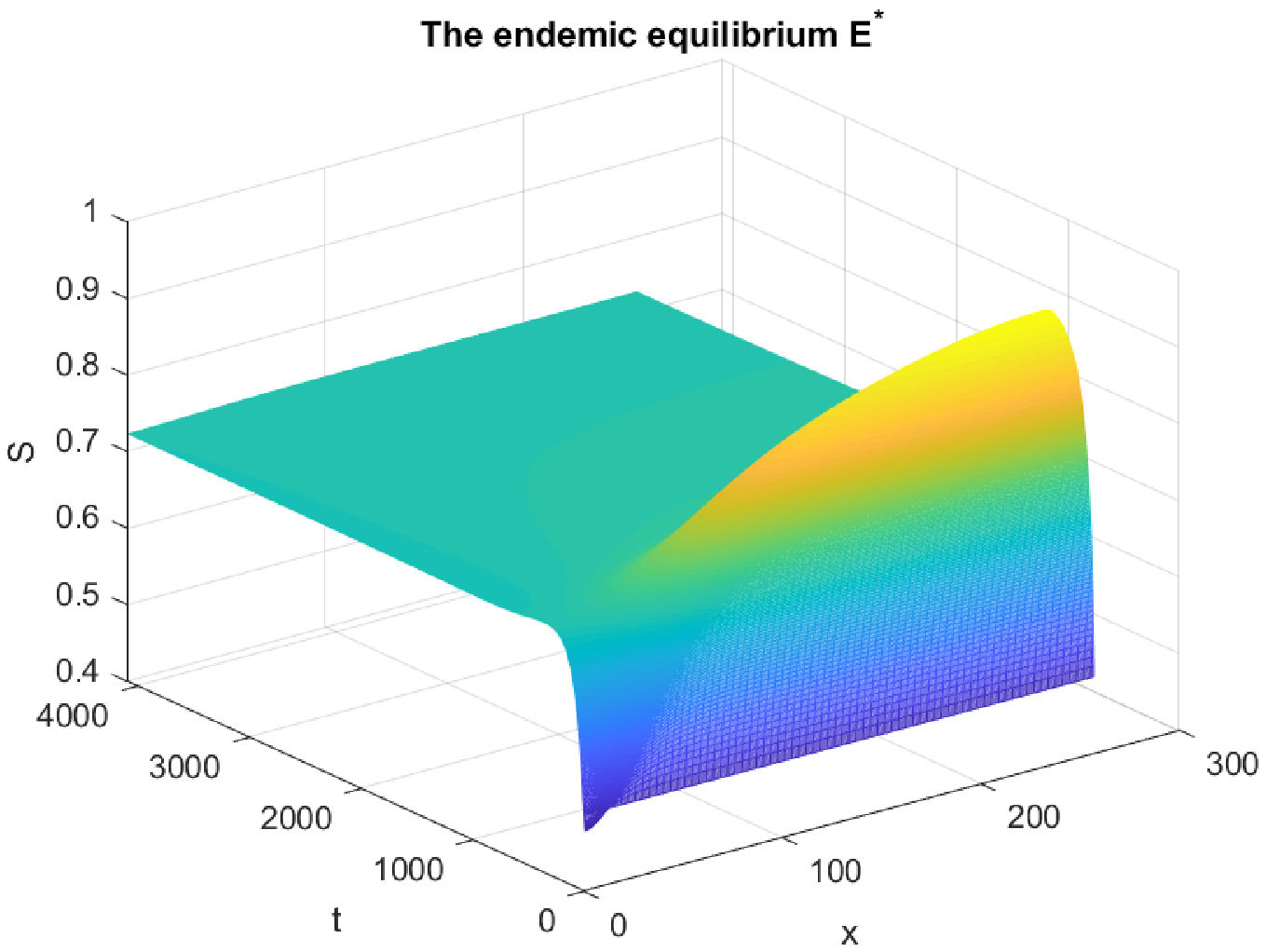}
       \includegraphics[width=0.4\linewidth]{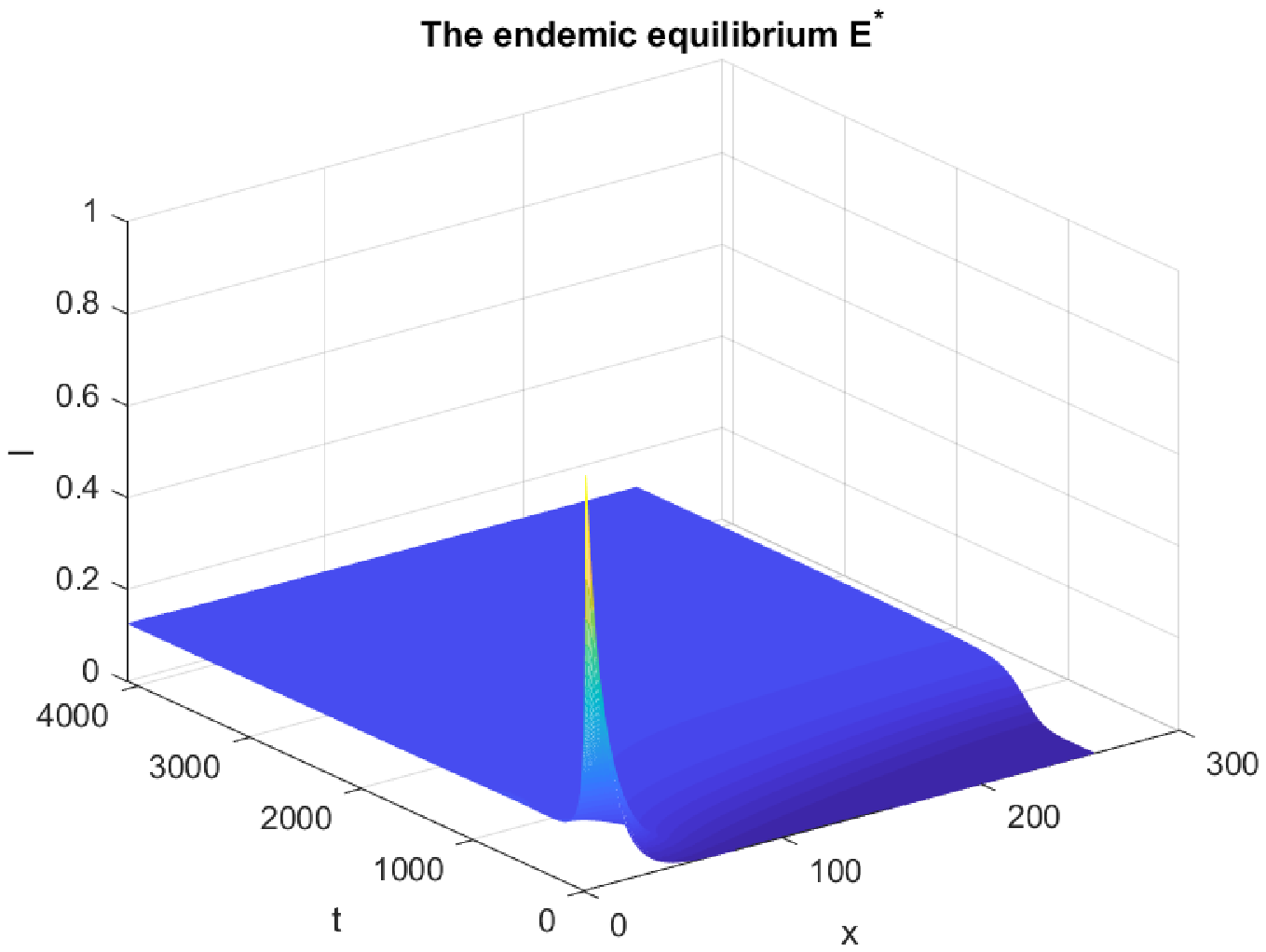}
       \includegraphics[width=0.4\linewidth]{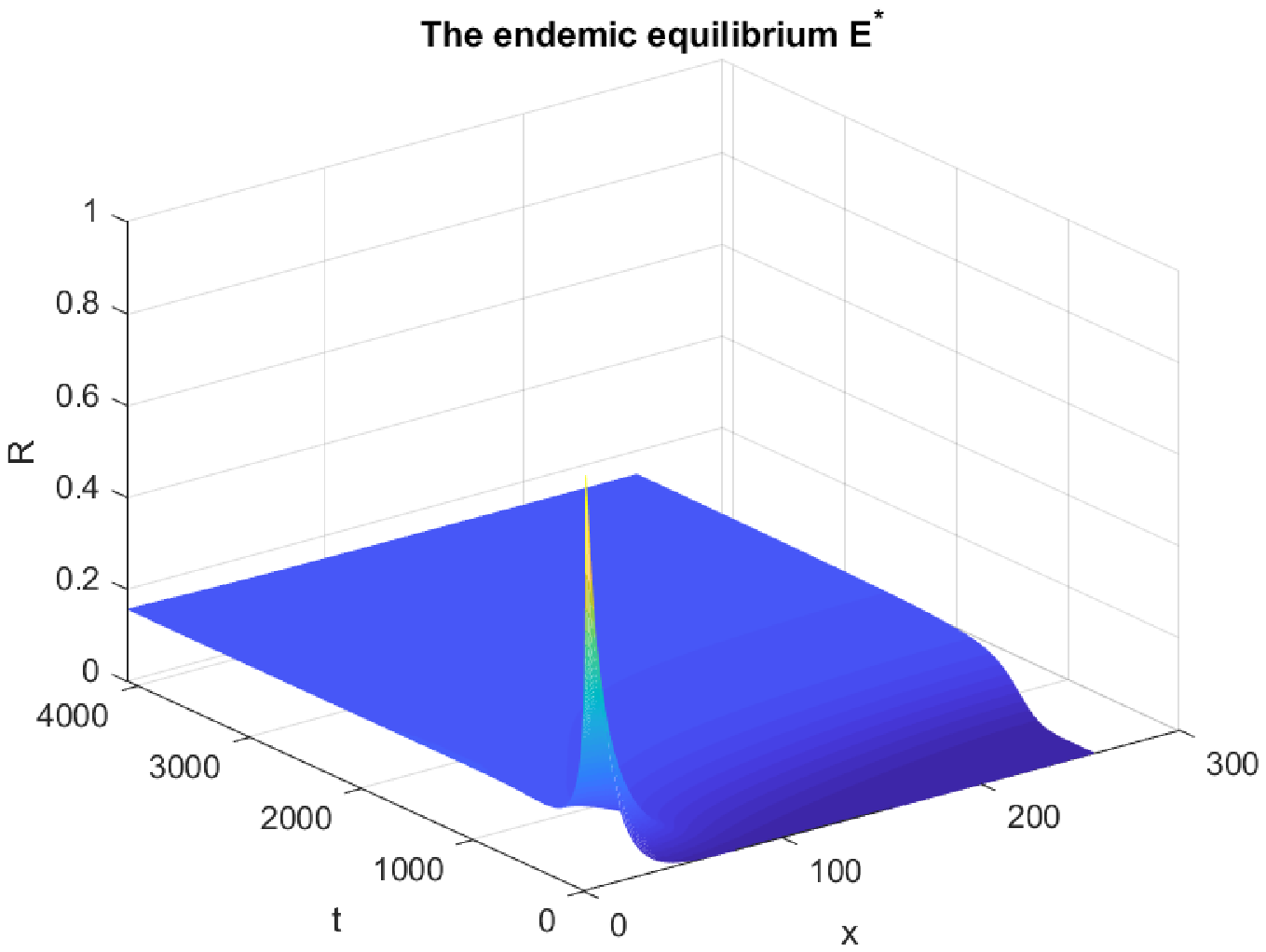}
      \caption{The endemic equilibrium $E^*$ is globally asymptotically stables.}
      \label{3}
    \end{center}
    \end{figure}

Next, we choose $\beta=0.6217$ and keep the other parameter values. The reproduction number $R_0=1.3816>1$ and the endemic equilibrium $E^*=(0.7238,0.1227,0.1534)$ by Lemma 3.3. It also follows from Theorem 3.4 that $E^*$ is globally asymptotically stable. This means that the disease persists and becomes chronic. Fig.~\ref{3} validates the above analysis.

\section{Conclusion}
In this paper, a discrete counterpart of time-fractional reaction-diffusion epidemic model with the generalized incidence rate is considered by second order central difference scheme and L1 nonstandard finite difference scheme. It should be noted that the  nonstandard difference scheme can not only keep the positive and boundedness of system (6), but also system (6) is a discrete system with time delay, and the time delay can not start from the initial value $t_0=0$. Furthermore, it can be found that the proposed system (6) can be simplified to the discrete system by forward difference scheme when the fractional-order $\alpha=1$. Meanwhile, the global properties of the discrete system (6) are analyzed and the results obtained are consistent with those of the corresponding continuous system. This shows that the discrete system composed of L1 nonstandard finite difference scheme and second-order central difference scheme can keep the property of continuous system well, including the positivity, boundedness and global stability of two equilibrium points, in which Lyapunov functions with memorization are firstly investigated to study the stability of discrete systems (6). Finally, the numerical simulations demonstrate the efficiency and applicability of the proposed system.

Furthermore, Lyapunov function constructed in this paper is consistent with the integer-order system as $\alpha\to 1^-$, which also ensures the $\alpha$-robustness of L1 scheme under the investigation of Lyapunov function. Meanwhile, Theorems 3.3 and 3.4 also guarantee the stability of integer-order discrete systems in \cite{2018Global} when $\alpha\to 1^-$. Therefore, in general, L1 scheme can not only keep the property of continuous system, but also is $\alpha$-robust in the Lypaunov function and the stability of the system. However, the establishment of systems from the perspective of discretization and other discrete schemes of Caputo fractional operator are not considered in this paper, which will be investigation in future.

\section*{Conflict of Interest}

The authors declare no conflicts of interest.

\section*{Acknowledgment}

\section*{References}

\bibliographystyle{elsarticle-num}
\bibliography{mybibfile}
\end{document}